\documentclass[11pt]{article}

\usepackage{amsmath,amsthm,verbatim,amssymb,amsfonts,amscd}
\usepackage {amsmath, amssymb, amscd, amsthm, epsfig, color}

\input xy
\xyoption {all}

\theoremstyle{plain}
\newtheorem{theorem}{Theorem}
\newtheorem{corollary}{Corollary}
\newtheorem{lemma}{Lemma}
\newtheorem{proposition}{Proposition}

\theoremstyle{definition}

\newcommand{\Z}{\mathbb{Z}}

\newcommand{\C}{\mathbb{C}}

\begin{document}
\title{Framed Rank $r$ Torsion-free Sheaves on $\C P^2$ and Representations of the Affine Lie Algebra $\widehat{gl(r)}$} 
\maketitle
\begin{abstract}
 We construct geometric realizations of the $r$-colored bosonic and fermionic Fock space on the equivariant cohomology of the moduli space of framed rank $r$ torsion-free sheaves on $\C P^2$.  Using these constructions, we realize geometrically all level one irreducible representations of the affine Lie algebra $\widehat{gl(r)}$.  The cyclic group $\Z_k$ acts naturally on the moduli space of sheaves, and the fixed-point components of this action are cyclic Nakajima quiver varieties .  We realize level k irreducible representations of $\widehat{gl(r)}$ on the equivariant cohomology of these quiver varieties.  \end{abstract}

\section{Introduction}
Let $M(r,n)$ be the moduli space of rank $r$ torsion free sheaves on $\C P^2 \subset \C P^2$, framed at 
$\C P^1_{\infty}$, with second Chern class $n$.  $M(r,n)$ is a partial compactification of the moduli space of $U(r)$ instantons on $\C^2$ with instanton number $n$.

For each $\vec{l} \in Z^r$, there is a one-dimensional torus $\C^*$ acting on $M(r,n)$.  Denote the space $M(r,n)$ with the torus action specified by $\vec{l}$ by $M_{\vec{l}}(r,n)$.  In addition to the one-dimensional torus action specified by $\vec{l}$, there is an $r$ dimensional torus which acts on $M(r,n)$ by changing the framing at $\C P^1_{\infty}$; all together this gives an action of an $r+1$ dimensional torus $T$ on the spaces $M_{\vec{l}}(r,n)$.  The $T$-equivariant cohomology of $M_{\vec{l}}(r,n)$ has a natural subspace 
$H^{mid}_T(M_{\vec{l}}(r,n),\C)$ whose dimension is equal to the Euler characteristic of $M_{\vec{l}}(r,n)$.  In this paper, we prove the following theorem:
\begin{theorem}	(i)	
\[
	\bigoplus_{\vec{l},n} H^{mid}_T(M_{\vec{l}}(r,n),\C)
\]
is a geometric realization of an irreducible module for an infinite dimensional Clifford algebra.
\\
(ii)	For any $\vec{l} \in \Z^r$, the space
\[
	\bigoplus_n H^{mid}_T(M_{\vec{l}}(r,n),\C)
\]
is a geometric realization of  an irreducible module for an infinite dimensional Heisenberg algebra.\\
(iii)	Let $Q \subset \Z^r$ is the sublattice $\{(l_0,\hdots,l_{r-1}) \mid \sum l_i = 0 \}$.

\[
	\bigoplus_{\vec{l} \in Q,n} H^{mid}_T(M_{\vec{l}}(r,n),\C)
\]
is a is geometric realization of the basic representation of $\widehat{gl(r)}$.
\end{theorem}

The above representations are constructed by exhibiting explicit correspondences inside products of the spaces $M(r,n)$.  In order to prove that these correspondences satisfy the correct relations, we make use of equivariant localization, which relates the equivariant cohomology of $M(r,n)$ to the equivariant cohomology of various submanifolds $Y \subset M(r,n)$ which are stable under the action of the torus $T$.  If we take $Y = M(r,n)^{T'}$, where $T'$ is the $r$-dimensional subtorus which changes the framing, then the connected components of $Y$ are isomorphic to products of Hilbert schemes ${\C^2}^{[n_0]} \times \hdots \times {\C^2}^{[n_{r-1}]}$.  We prove the Heisenberg commutation relations by pulling back our operators to operators on the equivariant cohomology of Hilbert schemes.  This part of the construction is thus a natural extension of the work of Nakajima [Na] and Grojnowski [Gr], as modified by Vasserot [Vas], who constructed geometric realizations of Heisenberg algebra representations using the geometry of Hilbert schemes.

Representations of Heisenberg and Clifford algebras are very closely related; in fact, starting from a repesentation of a Clifford algebra (a "fermionic" object), one can construct representations of a Heisenberg algebra (a "bosonic" object), and vice-versa.  The translation between the language of bosonic and fermionic operators, which was initially discovered by physicists, is known in the mathematics literature as the "boson-fermion correspondence."  Our geometric Heisenberg and Clifford algebra constructions provide a geometric interpretation of this correspondence.  This extends results of Savage [Sav], who studies the construction of a Heisenberg algebra representation  on the cohomology of the Hilbert Schemes ${\mathbb{C}^2}^{[n]}$, and relates this construction to a geometric realization of level one representations of the Lie algebra $sl(\infty)$.

In order to construct higher level representations of $\widehat{gl(r)}$, we add an action of the cyclic group $\Z_k$ to the spaces $M_{\vec{l}}(r,n)$.  The connected components of the fixed point set $M_{\vec{l}}(r,n)^{\Z_k}$ are Nakajima quiver varieties whose underlying graph is the 
Dynkin diagram of the affine Lie algebra $\widehat{sl(k)}$.  As vector spaces, there is a natural isomophism
\[
	\bigoplus_{\vec{l},n} H^{mid}_T(M_{\vec{l}}(r,n)^{\Z_k},\C) \simeq 
	\bigoplus_{\vec{l},n} H^{mid}_T(M_{\vec{l}}(r,n),\C).
\]
However, inside products of the spaces $M_{\vec{l}}(r,n)^{\Z_k}$ we can define correpondences on which make
\[
	\bigoplus_{\vec{l},n} H^{mid}_T(M_{\vec{l}}(r,n)^{\Z_k},\C)
\]
into a representation of a different infinite dimensional Heisenberg algebra.  Passing from a representation of this Heisenberg algebra to a representation of the Lie algebra $\widehat{gl(r)}$, we obtain the following theorem.
\begin{theorem}\label{higher level}
\[
	\bigoplus_{\vec{l},n} H^{mid}_T(M_{\vec{l}}(r,n)^{\Z_k},\C)
\]
is a geometric realization of a direct sum of irreducible level $k$ representation of $\widehat{gl(r)}$.
\end{theorem}

The connection between the representation theory of affine Lie algebras and instanton geometry was discovered by H. Nakajima in the remarkable work [Na94] (see also [Na98]). In [Na94], Nakajima constructed representations of infinite dimensional Lie algebras on the homology of quiver varieties, which are generalizations of $U(r)$-instanton moduli spaces.  In particular, Nakajima constructs level $r$ representations of an affine Lie algebra $\widehat{g}$ on the homology of moduli spaces of $U(r)$ instantons on $\widetilde{\C^2/\Gamma}$, where $\Gamma \subset SL(2,\C)$ is a finite subgroup and $\Gamma$ and $\widehat{g}$ are related by the McKay correspondence. In this construction the finite subgroup $\Gamma$ determines the Lie algebra being represented and and the gauge group $U(r)$ determines the level of the representation.

It seems equally natural, however, to expect algebraic objects associated to $G$ (rather than to the finite subgroup $\Gamma$) to be related to the topology of moduli spaces of $G$-instantons on $\widetilde{\C^2/\Gamma}$.  When $\Gamma = \Z_k$ is a cyclic group, a conjecture of I. Frenkel says that the homology of the moduli space of $G$-instantons on $\widetilde{\C^2/\Z_k}$ should carry level $k$ representations of the affine Lie algebra $\widehat{g}$ associated to $G$.  This second construction of affine Lie algebra representations is different from the Nakajima construction, in that the gauge group $G$ determines the algebra being represented, and the finite subgroup $\Z_k$ determines the level of the representation.\\
To summarize, fix a simply-laced affine Lie algebra $\widehat{g}$, corresponding to both a Lie group $G$ and a finite subgroup $\Gamma$.  There are two different constructions of representations of the Lie algebra $\widehat{g}$  on the homology of instanton moduli spaces:

(a) $\widehat{g}$ acts on the homology of moduli spaces of $U(k)$-instantons on $\widetilde{\C^2/\Gamma}$.  The level of the representation is determined by the gauge group $U(k)$; and\\

(b) $\widehat{g}$ acts on the homology of moduli spaces of $G$-instantons on $\widetilde{\C^2/\Z_k}$. The level of the representation is determined by the finite subgroup $\Z_k$.\\

The construction (a) is contained in [Na94], while this paper gives the construction (b) when $G=U(r)$ is of type A.  In this case the Lie algebra being represented is the affine Lie algebra $\widehat{gl(r)}$, and all of the moduli spaces involved have descriptions as Nakajima quiver varieties.  When $G \neq U(r)$, the construction (b) is still conjectural.  If $G \neq U(r)$ and $\Gamma \neq \Z_k$, it is an interesting question to determine what sort of exotic representations can be realized using the corresponding instanton moduli spaces.

We note that on the equivariant cohomology of the quiver varieties that Nakajima used to construct level $r$ representations of $\widehat{sl(k)}$, we construct level $k$ representations of $\widehat{gl(r)}$.  Thus, we have a level-rank duality in Nakajima quiver varieties of affine type $\widehat{A}$. In [Fr], Frenkel studied an algebraic level-rank duality in the representation theory of $\widehat{gl(r)}$.  We expect the contructions (a) and (b) above to offer geometric interpretations of algebraic level-rank duality.

Equivariant cohomology and localization have been used before in order to study the topology of the moduli space $M(r,n)$.  Of particular note are the papers [NY1], [NY2], which contain much of the information used in the constructions of this paper.  The goal of [NY1], [NY2], however,  which is to study instanton counting, does not require a geometric realization of representations. It should be interesting to relate instanton counting on surfaces to geometric realizations of representations on instanton moduli spaces.

\subsection{Acknowledgements}
The author is grateful to his advisor I. Frenkel for his guidance, generosity and encouragement during the duration of this project.  In addition, the author would like to thank Alina Marian, Josh Sussan, Manish Patnaik, Alistair Savage and Kevin Wortman for helpful conversations.

\section{Algebraic Preliminaries}
\subsection{Partitions and Symmetric Functions}
Let $\lambda = (\lambda_0 \geq \hdots \geq \lambda_m >0)$ be a partition of $n$, which we write using the notation $\lambda \vdash n$.\\  
We may associate to $\lambda$ its Young diagram $Y_{\lambda}$, which we view as a subset of the first quadrant of $\Z^2$, as in [NY1].  Given a positive integer $k$, we color the diagram $Y_{\lambda}$ with $k$ colors by coloring the node $(i,j)$ the color $j-i$ (mod $k$).\\
We say that $\lambda$ is $k$-regular if each color $0,\hdots, k-1$ occurs with the same multiplicity in the $k$-coloring of $Y_{\lambda}$.  Note that if $\lambda \vdash n$ is $k$-regular, then $n$ is divisible by $k$.\\
Let $Sym$ be the $\C$-vector space of symmetric functions.  Of the many important bases of this space, will have occasion to use the following  (see [Mac] for definitions):\\
the monomial symmetric functions $m_{\lambda}$\\
the power sum symmetric functions $p_{\lambda}$\\
the Schur functions $s_{\lambda}$\\
the elementary symmetric functions $e_{\lambda}$\\
the homogeneous symmetric functions $h_{\lambda}$.\\
$Sym \simeq \C[p_1,p_2,\hdots]$ is a polynomial algebra in the power-sum symmetric functions $\{p_n\}_{n >0}$.\\  
Denote by $Sym_{k-reg} \subset Sym$ the subspace spanned by the Schur functions $s_{\lambda}$ for   $k$-regular partitions $\lambda$.  

\subsection{The Clifford Algebra}
\noindent
Let $Cl$ be Clifford Algebra generated by $\psi(k),\psi^*(k)$, $k \in \Z$, and an element $c$, with anti-commutation relations
\[
\{\psi(k),\psi(l)\}=\{\psi^*(k),\psi^*(l)\}=0, \hskip.5cm
\{\psi(k),\psi^*(l)\}=\delta_{kl} c
\]

Define the spin module $\mathcal{F}$ to be the unique irreducible Clifford Module which admits a vector $\nu_0$ such that
 \[
 c \nu_0 = \nu_0
\]
\[
 \psi(k) \nu_0 = 0  \hskip.5cm \forall k \leq 0\]
\[\psi^*(k) \nu_0 = 0 \hskip.5cm  \forall k > 0\]
The Spin Module $\mathcal{F}$ is also known as the Fermionic Fock space, and this space has a nice realization in terms of semi-infinite monomials.  A semi-infinite monomial is an infinite expression of the form
\[
	i_0 \wedge i_2\wedge i_3 \wedge  \hdots
\]
where $i_0 >i_1>i_2>\hdots $ are integers and $i_n = i_{n-1} - 1$ for $n \gg 0$.\\
For any semi-infinite monomial $i_0\wedge i_1\wedge i_2 \wedge \hdots$, there exists
$k \in \Z$ such that for $n \gg 0$, $i_n = -n + k$, and we will refer to this $k$ as the charge of the semi-infinite wedge $i_0\wedge i_1\wedge i_2 \wedge \hdots $.  Put another way, the charge of $i_0\wedge i_1\wedge i_2 \wedge \hdots$ is the integer $k$ such that $i_0\wedge i_1\wedge i_2 \wedge \hdots$ differs from 
$k \wedge k-1 \wedge k-2 \wedge \hdots$ at only finitely many places.\\
Let $\mathcal{F}(m)$ be the $\C$-vector space spanned by all semi-infinite monomials of charge $m$, and let
\[
	\mathcal{F} = \bigoplus_{m} \mathcal{F}(m)
\]
The action of $\psi(k), \psi^*(k)$ on $\mathcal{F}$ is defined by wedging and contracting operators:
\[
	\psi(k)(i_0\wedge i_1 \wedge \hdots) = 
	 \left \{
   \begin{array}{ll}
      (-1)^s i_0\wedge \hdots \wedge i_{s-1} \wedge k \wedge i_s \hdots & \mbox{$i_{s-1} > k > i_s$}\\
     0 & \mbox{$k = i_s$ for some $s$}
    \end{array}
   \right.
\]
\[
\psi^*(k)(i_0\wedge i_1\wedge \hdots) =
 \left \{
   \begin{array}{ll}
    (-1)^s i_0\wedge \hdots \wedge {i_{s-1}} \wedge i_{s+1}\hdots& \mbox{$ k= i_s$}\\
     0 & \mbox{$k \neq i_s$ for all $s$}
    \end{array}
   \right.
\]
If we define an inner product on $\mathcal{F}$ by declaring the semi-infinite monomials to be an orthonormal basis, then $\psi(k)$ and $\psi^*(k)$ are adjoint operators.  Note that 
\[
	\psi(k) : \mathcal{F}(m) \longrightarrow \mathcal{F}(m+1)
\]
raises charge by one while 
\[
	\psi^*(k) : \mathcal{F}(m) \longrightarrow \mathcal{F}(m-1)
\]
lowers charge by one.  \\
There is also an $r$-tuple version of the Clifford Algebra, denoted by $Cl^r$; it is  generated by $\psi_i(k),\psi^*_i(k)$, $k \in \Z$, $i = 0,\hdots, r-1$ and an element $c$, with anti-commutation relations
\[
\{\psi_i(k),\psi_j(l)\}=\{\psi_i^*(k),\psi_j^*(l)\}=0, \hskip.5cm
\{\psi_i(k),\psi_j^*(l)\}=\delta_{ij}\delta_{kl} c
\]
We define the $r$-colored fermionic Fock space $\mathcal{F}^r$ by taking the tensor product of  $r$ copies of the space $\mathcal{F}$,
\[
	\mathcal{F}^r = \mathcal{F}\otimes \hdots \otimes \mathcal{F}
\]
so that $Cl^r$ acts naturally on $\mathcal{F}^r$.  \\
\vskip 0.5cm
\subsection{The Heisenberg Algebra}
\noindent
The Heisenberg Algebra $\mathcal{H}$ is the infinite dimensional Lie algebra generated by $c$, $p(n)$, $n \in \Z$, with commutation relations
\[
	[p(n),p(m)] = n\delta_{n+m,0}c
\]
\[
	[p(n),c] = 0
\]
Let $\mathcal{B}(k)$ be the unique irreducible $\mathcal{H}$ module which admits a vector $\nu_0$ such that
\[
	c \nu_0 = \nu_0
\]
\[
	p(n) \nu_0 = 0  \hskip.5cm \forall n < 0
\]
\[
	p(0) \nu_0 = k \nu_0
\]
and let $\mathcal{B} = \bigoplus_{k \in \Z} \mathcal{B}(k)$.\\
$\mathcal{B}$ is known as the Bosonic Fock space, and this space has a nice realization in terms of symmetric functions.\\
\vskip0.5cm
\noindent
Let $\mathcal{B}_{alg} = \C[q,q^{-1},p_1,p_2,\hdots]$, where $p_i$ are the power-sum symmetric functions, and let 
$\mathcal{B}_{alg}(k) = q^k \C[p_1,p_2,\hdots] = q^kSym$.
Define the action of $\mathcal{H}$ on the space on $\mathcal{B}_{alg}(k)$ by
\[
	cf = f
\]
\[
	p(n) f = p_n f\hskip.5cm \forall n > 0
\]
\[
	p(n) f = \frac{\partial}{\partial p_n} f \hskip0.5cm \forall n < 0
\]
\[
	p(0) f = k f
\]
Putting together the actions for different $k$, we have an action of $\mathcal{H}$ on $\mathcal{B}_{alg}$.
We define an inner product on $\mathcal{B}_{alg}$ by declaring the elements $q^ms_{\lambda}$ to be an orthonormal basis.  With this inner product, the operators $p(n)$ and $p(-n)$ are adjoints.\\
There is also an $r$-colored version of the Heisenberg algebra, denoted by $\mathcal{H}^r$.  This is the Lie algebra generated by $c$, $p_i(n)$, $n \in \Z$, $i = 0,\hdots,r-1$with commutation relations
\[
	[p_i(n),p_j(m)] = n\delta_{n+m,0}\delta_{i,j}c
\]
\[
	[p_i(n),c] = 0
\]
If we take $r$ copies of $\mathcal{B}_{alg}$ and set $\mathcal{B}^r_{alg} = \mathcal{B}_{alg} \otimes \hdots \otimes \mathcal{B}_{alg}$, then we have a natural action of $\mathcal{H}^r$ on $\mathcal{B}_{alg}^r$.  We will refer to $\mathcal{B}_{alg}^r$ as the $r$-colored bosonic Fock space.

\subsection{The Boson-Fermion Correspondence}
\noindent
We may associate a partition $\lambda = (\lambda_0 \geq \hdots \geq \lambda_k)$ to a semi-infinite monomial $i_0 \wedge i_1 \wedge \hdots$ of charge $m$ by setting
\[
	\lambda_j = i_j - m + j
\]
The correspondence
\[
	i_0 \wedge i_1 \wedge \hdots \leftrightarrow \lambda
\]	
allows us to define an isometric vector space isomorphism
\[
	\phi: \mathcal{F} \longrightarrow \mathcal{B}_{alg}
\]
\[
	\phi(i_0 \wedge i_1 \wedge \hdots) = q^ms_{\lambda}
\]
for a semi-infinite monomial $i_0 \wedge i_1 \wedge \hdots$ of charge $m$.\\
We use this isomorphism to define an action of $\mathcal{H}$ on $\mathcal{F}$, and an action of $Cl$ on $\mathcal{B}$.  In order to do so, we define the operators $h(k)$, $e(k)$ as homogeneous components of the generating functions
\[
	exp(\sum_{n=1}^{\infty} \frac{z^n}{n}p(n)) = \sum_{k=1}^{\infty} h(k) z^{k}
\]
\[
	exp(\sum_{n=1}^{\infty} \frac{z^n}{n}p(-n)) = \sum_{k=1}^{\infty} h(-k) z^{-k}
\]
and their inverses
\[
	exp(-\sum_{n=1}^{\infty} \frac{z^n}{n}p(n))  =  \sum_{k=1}^{\infty} e(-k) z^{-k}
\]
\[
	exp(-\sum_{n=1}^{\infty} \frac{z^n}{n}p(n)) = \sum_{k=1}^{\infty} e(k)z^{k}
\]
The operators $h(k), e(k)$ are adjoint to the operators $h(-k), e(-k)$, respectively. As operators on the space of symmetric functions, $h(k)$, $k>0$ is multiplication by the homogeneous symmetric function $h_k$.  Similarly, $e(k)$, $k>0$ is multiplicaton by the elementary symmetric funcion $e_k$ (see [Macdonald].)\\
We also define the shift operator 
\[
	q: \mathcal{B}_{alg}(k) \longrightarrow \mathcal{B}_{alg}(k+1)
\]
\[
	q(f) = qf
\]
Note that $q$ and $q^{-1}$ are adjoint operators.
\begin{proposition} [Fr 81]
a)	As operators on $\mathcal{B}_{alg}$, the bosons can be written in terms of the fermions:
\[
	p(n) = \sum_{j \in \Z} \psi(j+n)\psi^*(j)
\]
if $n \neq 0$, and
\[
	p(0) = \sum_{j >0} \psi(j)\psi^*(j) - \sum_{j \leq 0}\psi^*(j)\psi(j)
\]
b)	As operators $\mathcal{B}(m) \longrightarrow \mathcal{B}(m\pm1)$, the fermions can be written in terms of the bosons and the shift operator:
\[
	\psi(k) = \sum_{n\in \Z} q h(n)e(n-m + k)
\]
\[
	\psi^*(k) = \sum_{n \in \Z} q^{-1} e(n)h(n+m+k)
\]
\end{proposition}
\noindent
There is an $r$-colored version of this correspondence, using the isomorphism
$\mathcal{F}^r \simeq \mathcal{B}^r$, and $r$ different shift operators $q_0, \hdots q_{r-1}$.
We define the operators $h_i(k)$ to be homogeneous components of the generating function
\[
	exp(\sum_{n=1}^{\infty} \frac{z^n}{n}p_i(n)) = \sum_{k=1}^{\infty} h_i(k) z^{k}
\]
and similarly for the operators $e_i(k)$.  In terms of symmetric functions, the operator $h_i(k)$ is multiplication by the $i$th coordinate homogeneous symmetric function
$1\otimes\hdots \otimes h_k \otimes \hdots 1 \in Sym^r$, and similarly for the $e_i(k)$.

\begin{proposition}
a)	As operators on $\mathcal{B}_{alg}^r$, the bosons can be contructed from the fermions:
\[
	p_i(n) = \sum_{k \in \Z} \psi_i(k)\psi_i^*(k+n)
\]
if $n \neq 0$, and
\[
	p_i(0) = \sum_{j >0} \psi_i(k)\psi_i^*(k) - \sum_{j \leq 0}\psi_i^*(k)\psi_i(k)
\]
b)	As operators $\mathcal{B}_{alg}^r(m) \longrightarrow \mathcal{B}^r(m\pm1)$, the fermions can be constructed from the bosons and the shift operators:
\[
	\psi_i(k) = \sum_{n\in \Z} q_i h_i(n)e_i(n-m + k)
\]
\[
	\psi_i^*(k) = \sum_{n \in \Z} q_i^{-1} e_i(n)h_i(n+m+k)
\]
\end{proposition}

\subsection{The Affine Lie Algebra $\widehat{gl(r)}$}
Let $gl(r)$ denote the Lie algebra of $r \times r$ matrices.  Let $e_{i,j}$, $i,j = 0,\hdots,r-1$ denote the matrix units, which form a basis of the vector space $gl(r)$. Let $h_i = e_{ii}$ be the diagonal matrix units.  The span of the $\{h_i\}$ (the Cartan subalgebra of diagonal matrices) will be denoted by $\mathfrak{h}$.  If $n^+$ (resp. $n^-$) are the strictly upper triangular (resp. strictly lower triangular) matrices, then we have a decomposition
\[
	gl(r) = n^+ \oplus \mathfrak{h} \oplus n^-
\]
The $\Z$-lattice spanned by the $h_i$, called the weight lattice, will be denoted by $P$.\\
The affine Lie algebra $\widehat{gl(r)}$ is the infinite dimensional vector space
\[
	\widetilde{gl(r)} = gl(r) \otimes \C[t,t^{-1}] \oplus \C c
\]
with Lie bracket
\[
	[x \otimes t^k,y \otimes t^l] = [x,y]\otimes t^{k+l} + k\delta_{k+l,0}\langle x,y \rangle c
\]
where $\langle x,y \rangle= tr(xy)$ is the invariant bilinear trace form, and $c$ is central.\\
Let
\[
	\widehat{gl(r)} = \widetilde{gl(r)} \rtimes d
\]
be the semi-direct product of $\widetilde{gl(r)}$ with the derivative $d = t\frac{d}{dt}$.  The subalgebra
$\overline{\mathfrak{h}} = \mathfrak{h} \oplus \C c \oplus \C d$ is the Cartan subalgebra of $\widehat{gl(r)}$.  If we set 
\[
	\overline{n}^{\pm} = (n^{\pm} \oplus gl(r)) \otimes t^{\pm 1} \C[t^{\pm 1}]
\]
then we have a triangular decomposition
\[
	\widehat{gl(r)} = \overline{n}^+ \oplus \overline{\mathfrak{h}} \oplus \overline{n}^-.
\]
Let $\widehat{P} = P \oplus \Z c \oplus \Z d$ be the weight lattice of $\widehat{gl(r)}$.  The subset
\[
	\widehat{P}^{++} = \{a_{-1}d + a_0 h_0 + \hdots a_{r-1} h_{r-1} \mid a_{-1} \geq a_0 \geq \hdots \geq a_{r-1} \}
\]
is called the set of dominant weights.

\subsection{Highest Weight representations of $\widehat{gl(r)}$}
\noindent
For $\lambda \in \widehat{P}^{++}$, the irreducible highest weight representation $V(\lambda)$ is by definition the unique irreducible $\widehat{gl(r)}$ with a vector $\nu_0$ satisfying
\[
	n^+ \nu_0 = 0
\]
and
\[
	h\nu_0 = \langle \lambda, h \rangle \nu_0, \hskip0.4cm h \in \overline{\mathfrak{h}}.
\]
The integer $m = \langle \lambda,c \rangle$ is called the $level$ of the representation $V(\lambda)$.\\
\vskip0.5cm
One constructs highest weight representations of $\widehat{gl(r)}$ (see [Fr-K], [Se], [K-K-L-W]) from the bosonic or fermionic Fock space constructed above by using vertex operators to extend the action of the Heisenberg or Clifford algebra to an action of the entire affine Lie algebra.  For example, in order to construct the level one representations of $\widehat{gl(r)}$ inside the fermionic Fock space
$\mathcal{F}^r$, we introduce the normal ordering
\[
	:\psi_i(k)\psi_j^*(l): =  
	 \left \{
   \begin{array}{ll}
   \psi_i(k)\psi^*_j(l) & \mbox{if $j>0$}\\
    -\psi^*_j(l)\psi_i(k) & \mbox{if $j \leq 0$}
    \end{array}
   \right.
\]
We can then define an action of $\widehat{gl(r)}$ on $\mathcal{F}^r$ by setting
\[
	e_{i,j} \otimes t^k \mapsto \sum_{n \in \Z} :\psi_i(n+k)\psi_j^*(n):
\]
In particular, the $r$-dimensional Heisenberg algebra (called the "homogeneous Heisenberg subalgebra")
\[
	\mathcal{H}^r = \mathfrak{h} \otimes \C[t,t^{-1}] \oplus \C c
\]
acts on $\mathcal{F}^r$ as in the Boson-Fermion correspondence above
\[
	h_i \otimes t^k \mapsto \sum_{n \in \Z} :\psi_i(n+k)\psi_i^*(n):
\]
The spaces 
\[
	\mathcal{F}^r(m) = \sum_{m_0+ \hdots +m_{r-1} = m} \mathcal{F}(m_0) \otimes \hdots \otimes
	\mathcal{F}(m_{r-1})
\]
are all irreducible level one representations of $\widehat{gl(r)}$, and all of the irreducible level one representations of $\widehat{gl(r)}$ are realized as $\mathcal{F}^r(m)$ for some $m$.\\
Alternatively, one can start from the tensor product of the irreducible Heisenberg algebra representation $Sym^r$ and the lattice group algebra $\C[\Z^r]$, and construct the bosonic Fock space
\[
	\mathcal{B}^r = Sym^r \otimes \C[\Z^r]
\]
The operators given by the action of the Heisenberg algebra and translation in the lattice can be put together (using vertex operators) to define an action of $\widehat{gl(r)}$ on the space $\mathcal{B}^r$, which decomposes into irreducible level one representations.  The isomorphism between these two constructions (one on $\mathcal{F}^r$ and one on $\mathcal{B}^r$) essentially follows from the boson-fermion correspondence.\\

In order to get representations of level $k$ for $k > 1$, we consider the subalgebra
\[
	\widehat{gl(r)}_k \subset \widehat{gl(r)}
\]
\[
	\widehat{gl(r)}_k = gl(r) \otimes \C[t^k,t^{-k}] \oplus \C c.
\]
Note that there is a Lie algebra isomorphism $\widehat{gl(r)}_k \simeq \widehat{gl(r)}$. The space $\mathcal{F}^r$, when restricted to the subalgebra $\widehat{gl(r)}_k$, decomposes into a direct sum of irreducible level $k$ representations, and all irreducible level $k$ representations are realized in this way. \\

\section{Geometry of Quiver Varieties}
\subsection{Construction of Quiver Varieties}
Let $V$ be a an $n$-dimensional vector space, let W be an $r$-dimensional vector space, and let\\ 
$A,B \in Hom(V,V)$, $i \in Hom(W,V)$, and $j \in Hom(V,W)$.
Corresponding to the $\widehat{A}_0$ Dynkin diagram we have the varieties
\[
	\mathbb{M}(r,n) =
	\{(A,B,i,j) \mid[A,B] + ij=0 
	,stability\}\]
	\[M(r,n)=\mathbb{M}(r,n)/GL(n,\C).
\]
Here "stability" means that we take only those quadruples $(A,B,i,j)$ satisfying the following condition:\\
	If $i(W) \subset V' \subset V$ for some $A,B$-stable subset $V'$, then $V' = V$.  \\
This condition guarantees that $GL(n,\C)$ acts freely on $\mathbb{M}(r,n)$ (see [Na 94]).\\
\vskip0.5cm
\noindent
$M(r,n)$ is isomorphic to he moduli space of rank $r$ framed torsion-free sheaves on $\C P^2$ with second Chern class eqaul to $n$.  In the special case $r = 1$, $M(1,n)$ is isomorphic to the Hilbert Scheme of $n$ points on $\C^2$, which  we will denote by ${\C^2}^{[n]}$. (See [Na 99], Ch. 3.)\\
\vskip0.5cm
\noindent
For a positive integer $k$, let $\Z_k\subset SL(2,\C)$ be the cyclic subgroup of order $k$, let $\mathcal{R}_i$, $i=0,\hdots,k-1$ be the irreducible represenations of $\Z_k$, and let $Q$ be the two-dimensional $\Z_k$ module defined by the inclusion into $SL(2,\C)$.  We also allow $k= \infty$, in which case we set $\Z_{\infty} = \C^*$, embedded in $SL(2,C)$ as the diagonal matrices.  A pair of endomporphisms $A,B \in Hom(V,V)$ can be considered as a point $(A,B) \in Hom(Q \otimes V,V)$.
If in addition $V$ and $W$ are  $\Z_k$-modules, then $\Z_k$ acts on 
$\mathbb{M}(r,n)$ and $M(r,n)$. \\
Define
\[
	\mathbb{M}(\vec{w},\vec{v}) = \{(A,B,i,j) \in   Hom_{\Z_k}(Q \otimes V,V) \times
	Hom_{\Z_k}(W,V) \times Hom_{\Z_l}(V,W)
	\mid[A,B] + ij=0 
	,stability\}
\]
\[
	M(\vec{w},\vec{v}) = \mathbb{M}(\vec{w},\vec{v}) /\prod_{i} GL(V_i)
\]
Here $\vec{w} = (w_0,\hdots w_{k-1}),\vec{v} = (v_0,\hdots,v_{k-1})$ are the decomposition vectors of $W= \oplus_i W_i \otimes \mathcal{R}_i$ and $V = \oplus_i V_i \otimes \mathcal{R}_i$ into irreducbible $\Z_k$-modules, and  $Hom_{\Z_k}(X,Y)$ denotes the $\Z_k$-invariant part of $Hom(X,Y).$  \\
\noindent
The spaces $M(\vec{w},\vec{v})$ will be called $\widehat{A_{k-1}}$ quiver varieties (or $A_{\infty}$ quiver varieties in the case $l = \infty$,) and their relationship to the original space $M(r,n)$ is given by the following proposition:
\begin{proposition}
For a $\Z_k$ module $W \simeq \C^r$, the fixed point components of the natural $\Z_k$ action on $M(r,n)$ are the $\widehat{A_{k-1}}$ quiver varieties (or $A_{\infty}$ quiver varieties in the case $k = \infty$):
\[
	M(r,n)^{\Z_k} = \coprod_{\sum v_i = n} M(\vec{w},\vec{v})
\]
\end{proposition}
\begin{proof}
	See [Na97].
\end{proof}
\noindent
In the case $l=\infty$, $dimW = 1$, the non-empty $A_\infty$ quiver varieties are all points, and these points are parametrized by partitions $\lambda = (\lambda_1 \geq \hdots \geq \lambda_n)$.  Explicitly, let $W \simeq \C$ be the $\C^*$ module $\rho(t) = t^l$, and let $1_l$ be the decomposition vector of $W$, so that 
\[
	(1_l)_n = 
	 \left \{
   \begin{array}{ll}
     1 & \mbox{$l=n$}\\
     0 & \mbox{otherwise}
    \end{array}
   \right.
\]
\begin{proposition}
The collection of $A_{\infty}$ quiver varieties $\{M(1_l,\vec{v})\}$ is in natural bijective correspondence with the set of partitions of charge $k$.
\begin{proof}
	See  [Fr-Sa].
\end{proof}
\end{proposition}

\subsection{Torus actions on $M(r,n)$}
Let $T' = (\C^*)^r \subset GL(r,\C)$ be a maximal torus.  Let $\vec{a} = diag(a_0,a_1,\hdots,a_{r-1})\in T'$, and for 
$\vec{l} = (l_0,\hdots,l_{r-1}) \in \Z^r$ let $b_{\vec{l}} = diag(t^{l_0},t^{l_1},\hdots,t^{l_{r-1}})$.  We define an action of $T = \C^* \times T$ on $M(r,n)$ via
\[
	t (A,B,i,j) = (tA,t^{-1}B,i\vec{a}^{-1}b_{\vec{l}}^{-1},b_{\vec{l}}\vec{a}j)
\]	 
We denote the space $M(r,n)$ with the above $T$-action by $M_{\vec{l}}(r,n)$

\begin{lemma}
The fixed point components $M_{\vec{l}}(r,n)^{T'}$ are products of Hilbert Schemes on $\C^2$:
\[
	\coprod_n M(r,n)^{T'} = \coprod _{\vec{n}=(n_0,\hdots,n_{r-1})} 
	M_{l_0}(1,n_0) \times \hdots \times M_{l_{r-1}}(1,n_{r-1}) \simeq
	\coprod _{\vec{n}=(n_0,\hdots,n_{r-1})} {\C^2}^{[n_0]}
	\times \hdots \times {\C^2}^{[n_{r-1}]}
\]
\begin{proof}
See [NY1]
\end{proof}
\end{lemma}
\begin{lemma}
The $T$ fixed points in $M_{\vec{l}}(r,n)$ are isolated and naturally identified with the set of multipartitions
$\vec{\lambda} = (\lambda_0,\hdots,\lambda_{r-1})$ of charge $\vec{l} = (l_0,\hdots,l_{r-1})$ and total size $n$.
\begin{proof}
See [NY1]
\end{proof}
\end{lemma}
\noindent
The character of the Tangent Bundle to $M(r,n)$ at $\vec{\lambda}$, which we denote by $T_{\vec{\lambda}}$ is given by (see [NY1])
\[
	T_{\vec{\lambda}} = 
	\sum_{\alpha,\beta = 0}^{r-1} N_{\alpha,\beta}
\]
\[
	N_{\alpha,\beta}=
	t^{l_{\beta}-l_{\alpha}}e_{\beta}e_{\alpha}^{-1}
	\times \bigl{(}
	\sum_{s \in \lambda_{\alpha}} t^{-h_{\beta,\alpha}(s)} +
	\sum_{s \in \lambda_{\beta}} t^{h_{\alpha,\beta}(s)} \bigr{)}
\]
In the above formula, $h_{\beta,\alpha}(s)$ denotes the relative hook-length of the square $s \in \vec{\lambda}$ (See [NY1])\\
\noindent
The Tangent bundle to $M(r,n)^{T'}$ at $\vec{\lambda}$, which we denote by $U_{\vec{\lambda}}$ is given by taking the $\alpha=\beta$ part of $T_{\vec{\lambda}}$,
\[
	U_{\vec{\lambda}} = 
	\sum_{\alpha = 0}^{r-1} \bigl{(}
	\sum_{s \in \lambda_{\alpha}} t^{-h(s)} + 
	\sum_{s \in \lambda_{\alpha}} t^{h(s)} \bigr{)}.
\]
Here $h(s)$ is the ordinary hook length of $s$.

\subsection{$\widehat{A_{k-1}}$ quiver varieties and the Hilbert Scheme $\widetilde{\C^2/\Z_k}^{[n]}$}
In the special case $\vec{w} = (1,0,\hdots,0)$ and $\vec{v} = (n,n,\hdots,n)$, the $\widehat{A}_{k-1}$ quiver variety $M(\vec{w},\vec{v})$ is diffeomorphic to the Hilbert scheme $\widetilde{\C^2/\Z_k}^{[n]}$, discussed in detail in [Kuz].  We will denote the quiver variety $M(\vec{w},\vec{v})$ for this particular choice of $\vec{w},\vec{v}$ by $M(\vec{1_0},\vec{n})$.\\
There is an inclusion [G-K],[It-Na]
\[
	\widetilde{\C^2/\Z_k} \hookrightarrow {\C^2}^{[k]},
\]
and our $\C^*$-action on ${\C^2}^{[k]}$ restricts to a $\C^*$-action on $\widetilde{\C^2/\Z_k}$.  This, in turn, induces $\C^*$ actions on all of the Hilbert schemes $\widetilde{\C^2/\Z_k}^{[n]}$.
\begin{lemma}
The fixed points $(\widetilde{\C^2/\Z_k}^{[n]})^{\C^*}$ are isolated, and in natural bijection with the set of $k$-tuples of partitions $\vec{\lambda} = (\lambda_0,\hdots,\lambda_{k-1})$ of total size $n$.
\begin{proof}
	See [Kuz], [Q-W].
\end{proof}
\end{lemma}
Of course, the quiver variety $M(\vec{1_0},\vec{n})$ has a $\C^*$ action induced from the inclusion
\[
	M(\vec{1_0},\vec{n}) \hookrightarrow {\C^2}^{[kn]}
\]
and the fixed points $M(\vec{1_0},\vec{n})$ are given by the $k-regular$ partitions of $kn$.\\
Let $\lambda \in M(\vec{1_0},\vec{n})^{\C^*}$ be a $k$-regular partition of $kn$.  The character of the tangent space to 
$M(\vec{1_0},\vec{n})$ at $\lambda$ is given by taking the $\Z_k$=invariant part of the character of the tangent space to ${\C^2}^{[kn]}$ at $\lambda$, and is thus given by
\[
	T_{\lambda} = \sum_{s \in \lambda} t^{h(s)}\delta^k_{h(s),0} + \sum_{s \in \lambda} t^{-h(s)}\delta^k_{h(s),0}
\]
where $\delta^k_{a,b}$ is equal to 1 if $a \equiv b$ modulo $k$ and equal to $0$ otherwise.
The above formula shows that the only hook-lengths which appear as weights in the above tangent space are those which are divisible by $k$.  Letting $h^k(s) = \frac{h(s)}{k}$, we can re-write the above character as
\[
	T_{\lambda} = \sum_{s \in \lambda^k} t^{kh^k(s)} + \sum_{s \in \lambda^k} t^{-kh^k(s)}
\]
where $\lambda^k$ consists of those boxes of the Young diagram of $\lambda$ whose hook-lengths are divisible by $k$.
Now let $\vec{\lambda} = (\lambda_0,\hdots, \lambda_{k-1}) \in (\widetilde{\C^2/\Z_k}^{[n]})^{\C^*}$.  The character of the tangent space to $\widetilde{\C^2/\Z_k}^{[n]}$ at $\vec{\lambda}$ (see [Q-W]) is given by
\[
	T_{\vec{\lambda}} = \sum_{\alpha = 0}^{k-1} \big{(} \sum_{s \in \lambda_{\alpha}} t^{kh(s)} + 
	\sum_{s \in \lambda_{\alpha}} t^{-kh(s)} \big{)}
\]

A $\C^*$-equivariant diffeomorphism [Kuz]
\[
	f: M(\vec{1_0},\vec{n}) \longrightarrow \widetilde{\C^2/\Z_k}^{[n]}
\]
takes fixed points to fixed points, giving  an explicit combinatorial bijection between the fixed point sets $M(\vec{1_0},\vec{n})^{\C^*}$ and $(\widetilde{\C^2/\Z_k}^{[n]})^{\C^*}$.  Moreover, $f$ induces a $\C^*$-module isomorphism 
\[
	f_{*}: T_{\lambda} \longrightarrow T_{f(\lambda)}
\]
so that the above characters are actually the same.
As a function from $\mathcal{P}_{reg}(kn)$  to $\mathcal{P}^k(n)$, this bijection is known as the $k-quotient$ (see [Mac], [Lei]), a notion that occurs naturally in the modular representation theory of the symmetric group.  This diffeomorphism and the associated bijection will be important for us later when we discuss representations of level $k > 1$.

\subsection{Equivariant Cohomology of Quiver Varieties}
Let $M$ be a quiver variety of complex dimension $2m$, and suppose that $T=(\C^*)^k$ acts on $M$ with isolated fixed points.  Let $B_T$ be the classifying space of $T$, and let $E_T$ be the universal bundle.  $T$ acts freely on the space $E_T$, and hence freely on the product $M \times E_T$.  The equivariant cohomology of $M$ is defined to be the ordinary cohomology of the quotient space 
$M \times_T E_T$.
\[
	H^*_T(M,\C) = H^*(M \times_T E_T ,\C)
\]
We will always use complex coefficients for the equivariant cohomology of $M$, which we will denote by $H^*_{T}(M)$. \\
$H^*_{T}(M)$ is a module over $R = H^*_{T}(pt)$.  Let $\mathcal{R}$ denote the field of fractions of $R$, and let $\mathcal{H}_{T}^*(M) = H^*_{T}(M(r,n)) \otimes_{R} \mathcal{R}$ be the localized equivariant cohomology of $M$.\\
All of the usual cohomological constructions carry over to the equivariant setting.  In particular, if $V$ is a $T$-equivariant vector bundle on $M$, we have equivariant Chern classes $c_k(V) \in H^{2k}_T(M)$.  If $V$ is an $n$-dimensional vector bundle, the top equivariant Chern class $c_n(V)$ is called the 
equivariant Euler class of $V$, and is denoted by $e(V)$.\\
We endow $\mathcal{H}_{T}^*(M)$ with an inner product given by
\[
	\langle , \rangle : \mathcal{H}_{T}^*(M) \times \mathcal{H}_{T}^*(M) \longrightarrow \mathcal{R}
\]
\[
\langle x,y \rangle = (-1)^{m} p_*(i_*)^{-1}(x \cup y)
\]
where $i$ is the inclusion 
\[
i: M^{T} \hookrightarrow M
\]
and $p$ is the unique map from $M^{T}$ to a point
\[
p: M^T \longrightarrow \{pt\}
\]
For a $T$-stable smooth subvariety $Y \subset M$, the normal bundle $N_Y$ to $Y$ in $M$ is $T$-equivariant vector bundle.  
If a one-parameter subgroup $\C^* \longrightarrow T$ acts trivially on $Y$, then this subgroup induces a splitting
\[
	N_Y = N_Y^+ \oplus N_Y^0 \oplus N_Y^- 
\]
where 
\[
	N_Y^+ = \bigoplus_{n > 0} N_{Y}(n)
\]
is the positive weight space of $\C^*$,
\[
	N_Y^- = \bigoplus_{n < 0} N_Y(n)
\]
is the negative weight space of $\C^*$
and $N_Y^0$ is the zero weight space.\\
If $T' \subset T$ is the subgroup of $T$ which acts trivially on $Y$, then by choosing  a generic one-parameter subgroup of $T'$ we can guarantee that $N_Y^0 = 0.$
In general we can choose a splitting such that the equivariant Euler classes of the bundles $N_Y^+, N_Y^-$ satisfy
\[
	e(N_Y^+)= (-1)^m e(N_Y^-);
\]
for our purposes we will need to choose such a splitting in three different cases, as explained in the following three examples:  
\\
\textbf{Example 1} Let $M = M(r,n)$ and Let $T= \C^* \times T'$ be the $r+1$-dimensional torus acting on $M(r,n)$, where $T' \simeq (\C^*)^r$ be the $r$-dimensional subtorus acting only on the framing. Let
\[
	Y = M(r,n)^{T'} = \coprod_{(n_0,\hdots,n_{r-1}} {\C^2}^{[n_0]} \times \hdots \times {\C^2}^{[n_{r-1}]}
\]
Then the Normal bundle $N$ to $Y$ splits as a direct sum
\[
	N = \bigoplus_{\alpha \neq \beta =1 }^{r-1} N_{\alpha,\beta}
\]
and we choose the one-paramter subgroup given by $(1,1) \times (1,t,t^2,\hdots,t^{r-1}) \in \C^* \times T' = T$ so that
\[
	N^+ = \bigoplus_{\alpha < \beta} N^{\alpha,\beta}
\]
\[
	N^- = \bigoplus_{\alpha > \beta} N^{\alpha,\beta}
\]
Then, looking at the character for the tangent bundle, we see that $N_y^0 = 0$ and that $e(N^{\alpha,\beta}) = (-1) e(N^{\beta,\alpha})$, which implies that
\[
	e(N^+)= (-1)^{(r-1)n} e(N^-).
\]
\textbf{Example 2}  Fix a $\Z_k$ module $W \simeq \C^r$.  $\Z_k$ acts on $M(r,n)$ via the inclusion $\Z_k \hookrightarrow SL(2,\C)$ and the framing action $W$.  This $\Z_k$ action commutes with the $T$ action, so that the tori $T$ and $T'$ act on $M(r,n)^{\Z_k}$.  Let $M$ be a connected component of $M(r,n)^{\Z_k}$, and let 
\[
	Y \subset M^{T'}
\]
be a connected component of $(M(r,n)^{\Z_k})^{T'}$ (so that we have taken a $\Z_k$-invariant piece of example 1.)  Choose the one parameter subgroup of example one.  The normal bundle $N_k$ to $Y \subset M$ is given by taking the $\Z_k$-invariant part of the bundle $N$ in example one.  Again, looking at the character of this bundle, we have 
\[
	e(N_k^+)= (-1)^me(N_k^-).
\]
where $2m$ is the complex codimension of $Y$ in $M$.\\
\textbf{Example 3}  Let $M = M(r,n)$, with the $T$ action of example 1, an let
\[
	Y = M(r,n)^T.
\]
In this case, the fixed points are isolated, and the normal bundle at $y$ is the full tangent bundle $T_y$.  We choose the one-parameter subgroup given by $(t^r,t^{-r}) \times (1,t,t^2,\hdots,t^{r-1}) \in T$, which also has isolated fixed points.  Then, for $y \in Y$, we have
\[
		e(T_y^+)= (-1)^{rn} e(T_y^-).
\]

\subsection{Localization and the Transport Map $\eta$}
We return now to the general case of a smooth, $T$-stable subvariety $Y \subset M$.  Let $i_Y:Y \longrightarrow M$ be the inclusion. For $x \in H_{T}^*(Y)$, define $\eta(x)$ by
\[
	\eta(x) ={i_Y}_*(x) \cup e(N^-)^{-1} = \sum_j {i_{Y_j}}_* (x) \cup e(N^-_j)^{-1}
\]
where $\{Y_j\}$ are the connected components of $Y$.\\
A priori, the image of $\eta$ lies in the localized equivariant cohomology of $M$, since we divide by the equivariant Euler class.  However, using the argument on [Na] section $6$ (see also [Vas]), we have the following lemma:
\begin{lemma}
If $x \in H_{T}^*(Y)$ then $\eta(x) \in H^*_T(M)$.
\begin{proof}
Repeat the argument of [Na] section 6.
\end{proof}
\end{lemma}
Thus, we have a well-defined map
\[
	\eta: H_{T}^*(Y) \longrightarrow	H_{T}^*(M)
\]
\vskip0.5cm
\noindent
Define $H^{mid}_T(M)$ to be the image of $H_T^0(M^T)$ under the map $\eta$.
$H^{mid}_T(M)$ is a subspace of the middle dimensional equivariant cohomology $H^{2m}_T(M)$, and the dimension of $H^{mid}_T(M)$ is equal to the number of fixed points $\#\{M^T\}$.  This number is equal to the dimension of the total ordinary cohomology $H^*(M)$, since $M$ has a Bialinicki-Birula decomposition with one cell for each fixed point.  Since all the cells in the Bialinicki-Birula decomposition are even-dimensional, the quiver variety $M$ has no odd dimensional homology; thus the dimension of $H^{mid}_T(M)$ is the same as the Euler characteristic of $M$.\\
\begin{lemma}
Let $Y \subset M$ be a $T$-stable smooth subvariety, and let
\[
	\eta: H_T^*(M^T) \longrightarrow H^*_T(M)
\]
\[
	\eta': H^*_T(Y^T) \longrightarrow H^*_T(Y)
\]
\[
	\eta'': H_T^*(Y) \longrightarrow H^*_T(M)
\]
be the corresponding maps.  Then
\[
	\eta''(H^{mid}_T(Y)) \subset H^{mid}_T(M)
\]
\begin{proof}
	Let $y \in Y^T$, and let $N, N'$ be the tangent bundles to  $M$ at $y$ and to $Y$ at $y$, respectively.
	Then $N = N' \oplus N''$, where $N''$ is the normal bundle to $Y$ in $M$ at $y$.  It follows that
\[
	\eta = \eta'' \eta'.
\]
\end{proof}	
\end{lemma}

\begin{proposition}
The restriction of $\langle , \rangle$ to $H^{mid}_{T} (M)$ (note: non-localized) is non-degenerate and $\C$ valued.
\begin{proof}
	By the localization theorem ([C-K]), the classes $\eta(1_{\lambda})$ for points $\lambda \in M^{T}$ form a basis of
	$H^{mid}_{T^*} (M)$.  So, we directly compute
\begin{equation*}
	 \langle \eta(1_{\lambda}),\eta(1_{\mu}) \rangle = (-1)^m p_*(i_*)^{-1}(\eta(1_{\lambda}) \cup \eta(1_{\mu}) )=
\end{equation*}
\begin{equation*}
	= (-1)^m p_*(i_*)^{-1}({i_{\lambda}}_*(1_{\lambda}) \cup {i_{\mu}}_*(1_{\mu}) \cup e(N_{\lambda}^-)^-	\cup e(N_{\mu}^-)^-)=
	 \delta_{\lambda,\mu} (-1)^m e(N_{\lambda})^{-2} e(T_{\lambda}) = \delta_{\lambda,\mu}.
\end{equation*}
Thus, the classes $\eta(1_{\lambda})$ form an orthonormal basis, and the bilinear form restricted to 
$H^{mid}_{T} (M)$ is $\C$-valued and non-degenerate.
\end{proof}
\end{proposition}
\noindent
We will denote this restricion by $\langle , \rangle$ as well.
\begin{corollary}
	\[
		\eta: H^{mid}_T(Y) \longrightarrow  H^{mid}_T(M)
	\]
	is an isometry.
\begin{proof}
	Let 
	\[
		\eta_1: H^{0}_T(Y^T) \longrightarrow  H^{mid}_T(Y)
	\]
	\[
		\eta_2: H^{0}_T(Y^T) \longrightarrow  H^{mid}_T(M)
	\]
	Then, by the computation in the above lemma, $\eta_1$ and $\eta_2$ are isometries.  But
	$\eta_2 = \eta \eta_1$, and $\eta_1$ is surjective.  Thus $\eta$ is an isometry.
\end{proof}
\end{corollary}

\subsection{Localization of Correspondences}
Let $M_1,M_2$ be quiver varieties with a $T$-action, and let $Y_1 \subset M_1$, $Y_2 \subset M_2$ be $T$-stable smooth subvarieties.  Let $Z \subset Y_1 \times Y_2$ be a $T$-stable correspondence, so that the fundemental class $[Y]$ defines a linear map
\[
	[Z] : H^*_T(Y_1) \longrightarrow H^*_T(Y_2)
\]
\[
	[Z](a) = {q_2}_*(q_1^*(a) \cup [Z]).
\]
Here $q_1,q_2$ are he projections from $Y_1 \times Y_2$ to $Y_1,Y_2$ respectively.
We extend $\eta$ to a map on correspondences by setting
\[
	\eta([Z]) = (i_1 \times i_2)_* ([Z]) \cup (e(N_1^+)^{-1} \otimes e(N_2^-)^{-1})
\]
so that $\eta([Z])$ defines a linear map
\[
	\eta([Z]) : H^*_T(M_1) \longrightarrow H^*_T(M_2)
\]
\[
	\eta([Z])(b) = {p_2}_*(p_1^*(b) \cup \eta([Z])).
\]
Here $p_1,p_2$ are the projections from $M_1\times M_2$ to $M_1,M_2$, respectively.
\begin{proposition}
\[
	\eta([Z])(\eta(x)) = \eta([Z](x))
\]
\begin{proof}
Let 
\[
	p_j: M_1 \times M_2 \longrightarrow M_j \hskip0.4cm j=1,2
\]
and
\[
	q_j: Y_1 \times Y_2 \longrightarrow Y_j \hskip0.4cm j=1,2
\]
be the projection maps, and let
\[
	i_j : Y_j \longrightarrow M_j \hskip0.4cm j=1,2
\]
be the inclusion.  Then\\
$
\eta([Z])( \eta(x))=\\
= {p_2}_*\big( p_1^*(\eta(x))\cup \eta([Z]) \big)=\\
={p_2}_*\big( p_1^*({i_1}_*(x) \cup e(N_1^-)^{-1}) \cup (i_1 \times i_2)_*([Z]) \cup e(N_1^+)^{-1} \cup e(N_2^-)^{-1} \big)=\\
={p_2}_*\big( p_1^*({i_1}_*(x)) \cup (i_1 \times i_2)_*([Z]) \cup e(N_1)^{-1} \cup e(N_2^-)^{-1} \big)=\\
={p_2}_*(i_1 \times i_2)_* \big((i_1 \times i_2)^* p_1^*({i_1}_*(x)) \cup [Z] \cup e(N_1)^{-1} 
\cup e(N_2^-)^{-1} \big)=\\
={i_2}_*{q_2}_* \big( q_1^*i_1^*({i_1}_*(x)) \cup [Z] \cup e(N_1)^{-1} 
\cup e(N_2^-)^{-1} \big)=\\
={i_2}_*{q_2}_* \big( q_1^*(x) \cup e(N_1) \cup [Z] \cup e(N_1)^{-1} 
\cup e(N_2^-)^{-1} \big)=\\
={i_2}_*{q_2}_* \big( q_1^*(x) \cup [Z] \big) \cup e(N_2^-)^{-1} =\\
=\eta([Z](x))
$

\end{proof}
\end{proposition}
\noindent

\begin{corollary}	If $[Z] : H^{mid}_T(Y_1) \longrightarrow H^{mid}_T(Y_2)$ then
\[
	\eta([Z]) : H^{mid}_T(M_1) \longrightarrow H^{mid}_T(M_2)
\]
\begin{proof}
If $x \in H^{mid}_T(Y_1)$, then $\eta(x) \in H^{mid}_T(M_1)$.  By assumption
$[Z](x) \in H^{mid}_T(Y_2)$, so that $\eta([Z](x)) \in H^{mid}_T(M_2)$.
The claim now follows from the above proposition.
\end{proof}
\end{corollary}
The above proposition and corollary will allow us to compute the (anti-)comutation relations of an operator $[X]$ and it's adjoint $[X]^*$ on  $\bigoplus_n H^*_T(M(r,n))$ by computing relations of the transported operators  $\eta^{-1}([X]),\eta^{-1}([X]^*)$ on $\bigoplus_n H^*_T(Y(n))$ for suitably chosen $T$-stable subvarieties $Y(n) \subset M(r,n)$.\\
In particular, we will consider the subvarieties of examples $1,2,$ and $3$.  The inclusions
\[
	M(r,n)^{T'} \hookrightarrow M(r,n)
\]
\[
	(M(r,n)^{\Z_k})^{T'} \hookrightarrow M(r,n)^{\Z_k}
\]
\[
	M(r,n)^T \hookrightarrow  M(r,n)
\]
of examples $1,2,$ and $3$ above give rise to distinct transport maps
\[
	\eta_1:  H^{mid}_T(M(r,n)^{T'}) \longrightarrow H^{mid}_T(M(r,n))
\]
\[
	\eta_k: H^{mid}_T(M(r,n)^{\Z_k})^{T'}) \longrightarrow H^{mid}_T(M(r,n)^{\Z_k})
\]
\[
	\eta: H^{mid}_T(M(r,n)^T) \longrightarrow H^{mid}_T(M(r,n))
\]
which will be used to compute the commutation relations of the Heisenberg and Clifford operators defined in the next section.

\section{Geometric Realization of Heisenberg and Clifford Operators}
In this section we consider the fundemental vector space
\[
	\mathcal{B}^r = \bigoplus_{n,\vec{l}} H^{mid}_{T}(M_{\vec{l}}(r,n)).
\]
and define correpondences which induce operators
\[
	P_i(n),\psi_i(n),\psi_i^*(n) :\mathcal{B}^r \longrightarrow \mathcal{B}^r.
\]
These operators will make $\mathcal{B}^r$ into a module over the Heisenberg algebra $\mathcal{H}^r$ and the Clifford algebra $Cl^r$. \\

\subsection{The Clifford operators $\psi_i(k), \psi_i^*(k)$}
For simplicity, we begin with the case $r=1$.\\
For $x = (A,B,i) \in \mathbb{M}(1_l,\vec{v}), y = (A',B',i') \in \mathbb{M}(1_l,\vec{u})$, we write 
\[
	x \twoheadrightarrow y
\]
\noindent
if there exists $S \subset V$ an $A,B$-stable subset of dimension $\vec{v}-\vec{u} $ such that
\[
	(A_{V/S},B_{V/S},i_{V/S}) = (A',B',i').
\]
Here  $\mathbb{M}(1_l,\vec{v})$ is the affine space used to define the $A_{\infty}$ quiver variety
$ M(1_l,\vec{v})$, and $A_{V/S},B_{V/S}$ are the endomorphisms of the quotient space $V/S$ induced from $A$ and $B$.\\
\textbf{Example}:	The Hecke correspondence $E_k$ used by Nakajima [Na 94] to define the action of Chevalley generator $e_k \in sl(\infty)$ on $\bigoplus_{\vec{v}} (H^*(M(1,\vec{v})$ can be constructed in this notation as
\[
E(k) = \{(x,y) \in \mathbb{M}(1_l,\vec{v}) \times \mathbb{M}(1_l,\vec{v} + 1_k) \mid 
y \twoheadrightarrow x\}
\]
Then, modding out by the $\oplus_k GL(V_k)$ actions in both factors defines a correspondence $e_k$.  These correspondences and their adjoints generate the action of $sl(\infty)$ on the cohomology of $A_{\infty}$ quiver varieties.  See [Sav].\\
For a $\Z$-graded vector space $V$ of dimension $\vec{v}$ and homogeneous maps $A,B \in Hom(V,V)$, let $A(V)$ and $B(V)$ denote the images of the linear maps $A$ and $B$.  Let
$a(\vec{v})$ and $b(\vec{v})$ denote the dimension vectors of the vector spaces  $A(V), B(V)$ respectively.  \\
\noindent
If $x = (A,B,i) \in \mathbb{M}(1_l,\vec{v})$, note that $A(x) := (A_{A(V)},B_{A(V)},Ai)$ defines a point in 
$\mathbb{M}(1_{l+1},a(\vec{v}))$.\\
Similarly, $B(x) := (A_{B(V)},B_{B(V)},Bi)$ defines a point in $\mathbb{M}(1_{l-1},b(\vec{v})).$
\\
\noindent
Given a pair of integers $k>l$ let $\vec{k}^l$ be the vector
\[
k^l(i)=
 \left \{
   \begin{array}{ll}
     1 & \mbox{$l < i < k$}\\
     0 & \mbox{otherwise}
    \end{array}
   \right.
\] 

\vskip.5cm
\noindent
Define $\alpha(k)_{l,\vec{v}} \subset
\mathbb{M}(1_l,\vec{v}) \times \mathbb{M}(1_{l+1},\vec{v}+\vec{k}^l) $ as follows:
\[
\alpha(k)_{l,\vec{v}} = \{(x,y) \mid y \twoheadrightarrow A(x)\}
\]
Modding out by the $GL(V_k)$ actions in both factors, $\alpha(k)_{l,\vec{v}}$ defines a correspondence, which we denote by $\alpha(k)_{l,\vec{v}}$
\[
\alpha(k)_{l,\vec{v}} \subset  
M(1_l,\vec{v}) \times M(1_{l+1},\vec{v} + \vec{k}^l)
\]
Let $\beta(k)_{l,\vec{v}}$ denote the adjoint correspondence obtained by swapping the factors \\
 $M(1_l,\vec{v})$ and $M(1_{l+1},\vec{v} + \vec{k}^l)$.\\
 \noindent
Similarly, if $k \leq l$, let $\vec{k}_l$ be the vector
\[
k_l(i)=
 \left \{
   \begin{array}{ll}
     1 & \mbox{$k-1 < i < l$}\\
     0 & \mbox{otherwise}
    \end{array}
   \right.
\] 
Define $\beta(k)_{l,\vec{v}} \subset 
\mathbb{M}(1_l,\vec{v}) \times \mathbb{M}(1_{l-1},\vec{v}+\vec{k}_l) $ as follows:
\[
\beta(k)_{l,\vec{v}} = \{(x,y) \mid y \twoheadrightarrow B(x)\}
\]
Modding out by the $GL(V_k)$ actions in both factors, $\beta(k)_{l,\vec{v}}$ defines a correspondence
\[
\beta(k)_{l,\vec{v}} \subset 
M(1_l,\vec{v}) \times M(1_{l-1},\vec{v} + \vec{k}_l).
\]
Let $\alpha(k)_{l,\vec{v}}$ denote the adjoint correspondence obtained by swapping the factors \\
 $M(1_l,\vec{v})$ and $M(1_{l-1},\vec{v} + \vec{k}_l)$.\\
\vskip.5cm
\noindent
For $k,l \in \Z$,$\vec{v} = (v_m)_{m\in \Z}$, let 
 \[
 n(k,l,\vec{v}) =  \left \{
   \begin{array}{ll}
     v_k & \mbox{$k>l$}\\
     v_k + l - k & \mbox{$k\leq l$}
    \end{array}
   \right.
\]
Define operators $\psi(k)$, $\psi^*(k)$,$k \in \Z$ by
 \[
 \psi(k) = \bigoplus_{l \in \Z,\vec{v}} (-1)^{n(k,l,\vec{v})} \eta([\alpha(k)_{l,\vec{v}}])
 \]
 \[
 \psi^*(k) = \bigoplus _{l \in \Z,\vec{v}} (-1)^{n(k,l,\vec{v})} \eta([\beta(k)_{l,\vec{v}}])
 \]
where
\[
	\eta: \bigoplus_{n,l} H^{mid}_{\C^*}(M_l(1,n)^{\C^*}) \longrightarrow 
	\bigoplus_{n,l} H^{mid}_{\C^*}(M_l(1,n))
\]
is the transport map, extended to a map on correspondences as described in the section on equavariant
cohomology of quiver varieties.\\
On $r$-component products of $A_{\infty}$ quiver varieties we have $r$ different correspondences $\alpha_i(k)_{l,\vec{v}},\beta_i(k)_{l,\vec{v}}$, $i=0,\hdots r-1$ modifying only the $i$th factor of the product.
Accordingly, we can define operators 
\[
	\psi_i(k), \psi_i^*(k): \mathcal{B}^r \longrightarrow \mathcal{B}^r
\]
 using the map
 \[
	\eta: \bigoplus_{n,l} H^{mid}_{T}(M_{\vec{l}}(r,n)^{T}) \longrightarrow 
	H^{mid}_{T}(M_{\vec{l}}(r,n))
\]
Note that, by construction, the operators $\psi_i(k)$ and $\psi_i^*(k)$ are adjoint to one another. \\

\subsection{The Heisenberg operators $P_i(n)$}
Define $Z^o \subset \coprod_{n,k} {\C^2}^{[n]} \times  {\C^2}^{[k]} \times  {\C^2}^{[n+k]}$
to be the variety of triples $(A,B,C)$ such that $A$ and $B$ have disjoint support, and there exists an exact sequence
\[
	0 \rightarrow A \rightarrow C \rightarrow B \rightarrow 0
\]
Let $Z$ be the closure of $Z^o$.  The fundemental class $[Z]$ defines a multiplication
\[
	[Z]: \bigoplus_{n} H^{mid}_{\C^*}({{\C}^2}^{[n]}) \otimes \bigoplus_n H^{mid}_{\C^*}({{\C}^2}^{[n]})
	\longrightarrow \bigoplus_n H^{mid}_{\C^*}({{\C}^2}^{[n]})
\]
which makes $\bigoplus_{n} H^{mid}_{\C^*}({{\C}^2}^{[n]})$ into a commutative algebra (see [Gr]).
We construct an $r$-colored version of this algebra by defining
\[
	Z^{o,r} \subset \coprod_{\vec{n},\vec{k}} \big{(} {\C^2}^{[n_0]} \times \hdots \times   {\C^2}^{[n_{r-1}]} 		\big{)} \times \big{(}
	{\C^2}^{[k_0]} \times \hdots \times   {\C^2}^{[k_{r-1}]} \big{)}
	 \times  \big{(} {\C^2}^{[n_0 + k_0]} \times \hdots \times   {\C^2}^{[n_{r-1} + k_{r-1}]} \big{)}
\]
to be the variety of triples  $(\vec{A},\vec{B},\vec{C})$,$\vec{A} = (A_0,\hdots,A_{r-1})$,
$\vec{B} = (B_0,\hdots,B_{r-1})$, $\vec{C} = (C_0,\hdots,C_{r-1})$ such that for each $i$:\\
1) $A_i$ and $B_i$ have disjoint support, and\\
2) there exists an exact sequence
\[
	0 \rightarrow A_i \rightarrow C_i \rightarrow B_i \rightarrow 0
\]
Let $Z^r$ be the closure of $Z^{o,r}$. If
\[
	\eta_1:  \bigoplus_{n} H^{mid}_{T}(M_{\vec{l}}(r,n)^{T'}) \longrightarrow 
	\bigoplus_{n} H^{mid}_T(M_{\vec{l}}(r,n))
\]
is the transport map, then $\eta_1([Z])$ makes $\bigoplus_{n} H^{mid}_T(M_{\vec{l}}(r,n))$ into a commuative algebra.\\
Define a new $\C^*$ action on $M_{\vec{l}}(r,n)$ by
\[
	t \diamond (A,B,i,j) = (tA,B,i,tj).
\]
This $\C^*$ action commutes with the $T$ action on $M_{\vec{l}}(r,n)$.  In particular, we have an induced action of $\C^*$ (also denoted by $\diamond$) on $M_{\vec{l}}(r,n)^{T'}$.  The fixed point components $M_{\vec{l}}(r,n)^{T' \times \C^*}$ (which are not in general isolated) are naturally enumerated by $r$-tuples of paritions of total size $n$ (see [Na 99]), and we denote the fixed point component corresponding to $\vec{\lambda}$ by $C_{\vec{\lambda}}$.\\
For $\mu \vdash n$, let $\vec{\mu}_i$
be the $r$-tuple of partitions $(0,\hdots,0,\mu,0,\hdots,0)$ which is equal to $\mu$ in the $ith$ place and  empty otherwise.  For example,  $\vec{1^n}_i = (0,\hdots,0,1^n,0,\hdots,0)$, 
$\vec{n}_i = (0,\hdots,n,\hdots,0).$\\
For $n>0$, define classes $p_i(n),e_i(n), h_i(n) \in H^{mid}_T(M(r,n)$ by
\[
	p_i(n) =  [closure\{z \in M_{\vec{l}}(r,n)^{T'} \mid \lim_{t \to 0} t\diamond z \in C_{\vec{n}_i} \}]
\]
\[
	e_i(n) =  [closure\{z \in M_{\vec{l}}(r,n)^{T'} \mid \lim_{t \to 0} t\diamond z \in C_{\vec{1^n}_i} \}]
\]
\[
	h_i(n) =  [closure\{z \in M_{\vec{l}}(r,n)^{T'} \mid \lim_{t \to 0} t\diamond z \hskip0.3cm exists\}]
\]
Define the classes $P_i(n),H_i(n),E_i(n)$ by
\[
	P_i(n) = \eta_1(p_i(n))
\]
\[
	E_i(n) = \eta_1(e_i(n))
\]
\[
	H_i(n) = \eta_1(h_i(n))
\]
Multiplication by $P_i(n)$ defines an operator, also denoted by $P_i(n)$.  Denote the adjoint of this operator by $P_i(-n)$.  The operators $P_i(n)$, $n \in \Z$ will be our Heisenberg operators.  The operators $E_i(n)$, $H_i(n)$ will be important for us later when we discuss the boson-fermion correspondence.

\section{The Proof of the Commutation and Anti-commutation Relations}
\subsection{The Clifford Algebra}
Let  $\nu_0 = 1 \in H^0_{T}(M(r,0).$
\begin{proposition}
 The operators $\psi_i(k),\psi^*_i(k)$ satisfy the following anti-commutation relations:
 \[
 	\psi_i(k) \nu_0 = 0  \hskip.5cm \forall k \leq 0,i=0,\hdots,r-1\]
\[\psi_i^*(k) \nu_0 = 0 \hskip.5cm  \forall k > 0,i=0,\hdots,r-1\]
\[\{\psi_i(k),\psi_j(l)\}=\{\psi_i^*(k),\psi_j^*(l)\}=0, \hskip.5cm
\{\psi_i(k),\psi_j^*(l)\}=\delta_{ij}\delta_{kl}
\]
\end{proposition}
\begin{proof}
Let 
\[
	\mathcal{F}^r = \bigoplus_{n,\vec{l}} H^0_T(M_{\vec{l}}(r,n)^T)
\]
and let $\eta$ be the isomorphism
\[
	\eta: \mathcal{F}^r \longrightarrow \mathcal{B}^r.
\]
We will consider the operators 
\[
	\eta^{-1}(\psi_i(k)), \eta^{-1}(\psi_i^*(k)): \mathcal{F}^r \longrightarrow \mathcal{F}^r
\]
and prove that they have the above commutation relations.  For convenience, we will drop the notation $\eta^{-1}$, and denote our transported operators by $\psi_i(k), \psi_i^*(k)$ as well.  In addition, since $\psi_i(k),\psi^*_i(k)$ and $\psi_j(k),\psi^*_j(k)$ for $i \neq j$ act on different coordinates in the product of $A_{\infty}$ quiver varieties, the only interesting case is the case $i=j$; thus we may consider the case $r$ = 1, and drop the subscripts $i,j$.\\
In order to prove this proposition for $r=1$, it will be convenient to identify the vector space $\mathcal{F}$ with the semi-infinite wedge space.  Given an $A_{\infty}$ quiver variety $M(1_l,\vec{v})$, 
we define subsets $C^+_{\vec{v},l},C^-_{\vec{v},l}, C_{\vec{v},l} \subset \Z$ by 
\[
	C^+_{\vec{v},l} = \{k>l \mid v_k \neq v_{k-1} \}
\]
\[
	C^-_{\vec{v},l} = \{k \leq l \mid v_k = v_{k-1} \}
\]
\[
	C_{\vec{v},l} = C^+_{\vec{v},l} \cup C^-_{\vec{v},l}
\]	
Arranging the elements of $C_{\vec{v},l}$ in decending order, \\
$C_{\vec{v},l} = \{i_0,i_1,\hdots\}$\\
we get a semi-infinite wedge
\[
	C_{\vec{v},l} \mapsto i_0 \wedge i_1 \wedge i_2 \wedge \hdots
\]
We define the charge of a semi infinite wedge $i_0 \wedge  i_1 \wedge i_2 \wedge \hdots$ to be the integer $m$ such that $i_n = m - n$ for $n$ sufficiently large; in this way the quiver variety $M(1_l,\vec{v})$ corresponds to a semi-infinite wedge of charge $l$.
Let $F_l$ denote the $\C$-span of the semi-infinite wedges of charge $l$.  We define a vector space isomorphism 
\[
	\mathcal{F}_l = \oplus_{\vec{v}} H^0_{\C^*}(M(1_l,\vec{v})) \longrightarrow F_l
\]
by mapping $1 \in H^0_{\C^*}(M(1_l,\vec{v})$ to the semi-infinite wedge corresponding to $M(1_l,\vec{v})$.
The following lemma, which is easy to check, relates to coordinate entries $v_k$ of $\vec{v}$ to the integers appearing in the corresponding semi-infinite monomial. 
\begin{lemma} If
$M(1_l,\vec{v})$ corresponds to the wedge $i_0 \wedge i_1 \wedge i_2 \wedge \hdots$
then the number of elements in the set $\{i_0,i_1,i_2,\hdots\}$ which are greater than $k$ is\\
\[
  \begin{array}{ll}
     v_k & \mbox{if $k>l$}\\
     v_k + l - k & \mbox{if $k \leq l$.}
    \end{array}
 \]   
\end{lemma}  
\noindent
The anti-commutation relations of the operators $\psi(k),\psi^*(k)$ will follow immediately from the following
\begin{lemma}
\[
	\psi(k)(i_0\wedge i_1 \wedge \hdots) = 
	 \left \{
   \begin{array}{ll}
      (-1)^s i_0\wedge \hdots \wedge i_{s-1} \wedge k \wedge i_s \hdots & \mbox{$i_{s-1} > k > i_s$}\\
     0 & \mbox{$k = i_s$ for some $s$}
    \end{array}
   \right.
\]
\[
\psi^*(k)(i_0\wedge i_1\wedge \hdots) =
 \left \{
   \begin{array}{ll}
    (-1)^s i_0\wedge \hdots \wedge {i_{s-1}} \wedge i_{s+1}\hdots& \mbox{$ k= i_s$}\\
     0 & \mbox{$k \neq i_s$ for all $s$}
    \end{array}
   \right.
\]

\end{lemma}
\begin{proof}
 Fix $l\in\Z$.  
 It suffices to show that the proposition holds for the operators $\psi(k)$, $k>l$ and $\psi^*(k)$, $k\leq l$.\\
 We begin with the operators $\psi(k)$, $k>l$.\\
Suppose $M(1_l,\vec{v})$ is nonempty and
 $M(1_l,\vec{v}) \leftrightarrow (i_0\wedge i_1\wedge \hdots \wedge i_n \wedge \hdots)$. \\
Then $M(1_{l+l}, \vec{v} + \vec{k}^l)$ is non-empty if and only if $k \notin \{i_0,i_1,i_2 \hdots\}$.  If $M(1_{l+1}, \vec{v} + \vec{k}^l)$ is non-empty, then the projections 
\[
p_1: \alpha(k) \longrightarrow
	M(1_l, \vec{v})
\]
\[
p_2:\alpha(k) \longrightarrow
M(1_{l+1}, \vec{v} + \vec{k}^l)
\]
induce identity maps on cohomology, since the varieties involved are all points.\\
It follows that
\[
[\alpha(k)](i_0\wedge i_1\wedge \hdots \wedge i_n \wedge \hdots) = 
i_0\wedge i_1\wedge \hdots \wedge i_s \wedge k \wedge i_{s+1} \wedge \hdots =
(-1)^{s+1} k \wedge i_0\wedge i_1\wedge \hdots \wedge i_n \wedge \hdots
\]
where $i_s > k \geq i_{s+1}$.
Therefore
\[
\psi(k)(i_0\wedge i_1\wedge \hdots \wedge i_n \wedge \hdots) = 
(-1)^{v_k}[\alpha(k)](i_0\wedge i_1\wedge \hdots \wedge i_n \wedge \hdots)
 \]
 But by the previous lemma $v_k$ is the number of elements in the set $\{i_0,i_1,\hdots\}$ which are greater than $k$, i.e.
 $v_k = s$, where $i_{s-1} > k > i_s$.\\
 \[
 \Longrightarrow
\psi(k)(i_0\wedge i_1 \wedge \hdots) = 
	 \left \{
   \begin{array}{ll}
      (-1)^s i_0\wedge \hdots \wedge i_{s-1} \wedge k \wedge i_s \hdots & \mbox{$i_{s-1} > k > i_s$}\\
     0 & \mbox{$k = i_s$ for some $s$}
    \end{array}
   \right.
\]
\vskip0.5cm
\noindent
Now suppose that $k \leq l$.  We will consider the operators $\psi^*(k)$.
We have that $M(1_{l+1},\vec{v} + \vec{k}_l)$ is non-empty if and only if $k \in \{i_0,i_1,i_2,\hdots\}$.
As in the case $\psi(k)$,$k>l$, we have
\[
\psi^*(k)(i_0\wedge i_1\wedge \hdots \wedge i_n \wedge \hdots) = 
(-1)^{v_k+l-k}[\beta(k)](i_0\wedge i_1\wedge \hdots \wedge i_n \wedge \hdots)=
 \sum_{j\in \Z} \delta_{i_j,k}(-1)^{v_k + l - k} i_0\wedge i_1\wedge \hdots \wedge \widehat{i_j} \wedge \hdots
\]
But $v_k + l - k$ is the number of elements in $\{i_0,i_1,i_2,\hdots\}$ which are greater than $k$
\[
\Longrightarrow 
\psi^*(k)(i_0\wedge i_1\wedge \hdots \wedge i_n \wedge \hdots) =
\sum_{j\in \Z} \delta_{i_j,k}(-1)^j i_0\wedge i_1\wedge \hdots \wedge \widehat{i_j} \wedge \hdots =
\]
\[
 \left \{
   \begin{array}{ll}
    (-1)^s i_0\wedge \hdots \wedge {i_{s-1}} \wedge i_{s+1}\hdots& \mbox{$ k= i_s$}\\
     0 & \mbox{$k \neq i_s$ for all $s$}
    \end{array}
   \right.
\] 
 \end{proof}
This completes the proof of the proposition. 
\end{proof}
Thus we obtain part (a) of the theorem stated in the introduction.
\begin{theorem}\label{Clifford} 
\[
	\bigoplus_{\vec{l},n} H^{mid}_T(M_{\vec{l}}(r,n),\C)
\]
is a geometric realization of an irreducible module for the Clifford algebra $Cl^r$.
\end{theorem}

\subsection{The Heisenberg Algebra}
\begin{proposition}
The operators $P_i(n)$, $n \in \Z$, $i = 0,\hdots r-1$ satisfy 
\[
	P_i(n) \nu_0 = 0, \hskip0.3cm n<0
\]
\[
	[P_i(n),P_j(m)] = n \delta_{i,j} \delta_{n+m,0} Id
\]
This makes $\mathcal{B}^r$ into a geometric realization of the $r$-colored bosonic Fock space.
\end{proposition}
\begin{proof}
Recall the isomorphism
\[
	\eta_1: \bigoplus_{n} H^{mid}_T(M_{\vec{l}}(r,n)^{T'}) \longrightarrow H^{mid}_T(M_{\vec{l}}(r,n))
\]
We will consider the operators $\eta_1^{-1}P(n)$, and show that they satisfy the same relations.  For simplicity we drop the notation $\eta_1^{-1}$ and denote the transported operators by $P_i(n)$.  
Since $P_i(n)$ and $P_j(m)$ act on different factors of $M_{\vec{l}}(r,n)^{T'}$, the only interesting case is $i=j$; thus we may consider the case $r=1$.\\
For a partition $\lambda \in ({\C^2}^{[n]})^{\C^*}$, we have a class 
$[\lambda] = \eta(1_{\lambda}) \in H^{mid}_{\C^*}({\C^2}^{[n]})$.  The classes $\{[\lambda]\}_{\lambda \vdash n}$ form an orthonormal basis of $H^{mid}_{\C^*}({\C^2}^{[n]})$, so we can construct an isometric vector space isomorphism
\[
	\phi:  Sym \longrightarrow \bigoplus_n H^{mid}_{\C^*} ({\C^2}^{[n]})
\]
\[
	\phi(s_{\lambda}) = [\lambda]
\]
by sending the Schur function $s_{\lambda}$ to the class  $[\lambda]$.  Then we have the following two lemmas.
\begin{lemma}
The map $\phi$ is an isomorphism of algebras.
\begin{proof}
See [Na],[Vas].
\end{proof}
\end{lemma}
\begin{lemma}
$\phi(m_{\lambda}) = [L_{\lambda}]$
\begin{proof}
See [Na],[Vas].
\end{proof}
\end{lemma}
In particular, the operator $P(n)$, $n>0$ is the image under $\phi$ of multiplication by the power-sum symmetric function $p_n$.  Since $\phi$ is an isometry, the adjoint operator $P(-n)$ corresponds to the differential operator $\partial/\partial p_n$.  The proposition then follows from
\[
	\partial/\partial p_n (1) = 0
\]
and
\[
	[p_n,\partial/\partial p_m] = n \delta_{n,m} Id
\]
\end{proof}
This gives part (b) of the theorem stated in the introduction.
\begin{theorem}\label{Heisenberg}	For any $\vec{l} \in \Z^r$,
\[
	\bigoplus_n H^{mid}_T(M_{\vec{l}}(r,n),\C)
\]
is a geometric realization of an irreducible module for the Heisenberg algebra $\mathcal{H}^r$.
\end{theorem}
Passing from representations of the Heisenberg algebra $\mathcal{H}^r$ to a representation of the affine Lie algebra $\widehat{gl(r)}$, we can identify the basic representation of $\widehat{gl(r)}$, as in  part (c) of the theorem stated in the introduction.  Let $Q \subset \Z^r$ be the sublattice whose entries sum to $0$.
\begin{theorem}\label{basic representation}
\[
	\bigoplus_{\vec{l} \in Q,n} H^{mid}_T(M_{\vec{l}}(r,n),\C)
\]
is a geometric realization of the basic representation of $\widehat{gl(r)}$.
\end{theorem}

\subsection{The boson-fermion Correspondence}
Now that we have constructed the fermionic operators $\psi_i(k), \psi_i^*(k)$ and the bosonic operators $P_i(n)$ on the common space
\[
	\mathcal{B}^r = \bigoplus_{n,\vec{l}} H^{mid}_{T}(M_{\vec{l}}(r,n))
\]
we can relate them to one another, giving a geometric realization of the boson-fermion correspondence.  The result in this section should be compared to the main result of [Sav], which considers the case $r=1$ and uses localization to relate the action of the Heisenberg algebra with an action of $sl(\infty)$.  We remark that the $sl(\infty)$ action considered in [Sav] can be constructed from the Clifford algebra action considered above.\\
\vskip0.5cm

\begin{theorem}
a)	As operators on $\mathcal{B}^r$, the bosons can be written in terms of the fermions:
\[
	P_i(n) = \sum_{k \in \Z} \psi_i(k)\psi_i^*(k+n)
\]
if $n \neq 0$, and
\[
	P_i(0) = \sum_{j >0} \psi_i(k)\psi_i^*(k) - \sum_{j \leq 0}\psi_i^*(k)\psi_i(k)
\]
b)	As operators $\mathcal{B}^r(m) \longrightarrow \mathcal{B}^r(m\pm1)$, the fermions can be written in terms of the bosons and the shift operators:
\[
	\psi_i(k) = \sum_{n\in \Z} Q H_i(n)E_i(n-m + k)
\]
\[
	\psi_i^*(k) = \sum_{n \in \Z} Q^{-1} E_i(n)H_i(n+m+k)
\]
\begin{proof}
Recall the algebraic version of the $r$-colored bosonic Fock space 
\[
	\mathcal{B}^r_{alg} = \mathcal{B}_{alg} \otimes \hdots \otimes \mathcal{B}_{alg},
\]
 where 
\[
	\mathcal{B}_{alg} = \C[q,q^{-1},p_1,p_2,\hdots].
\]
We denote the vector $q_0^{l_0}s_{\lambda_0} \otimes \hdots \otimes q^{l_{r-1}}s_{\lambda_{r-1}} \in \mathcal{B}^r_{alg}$
by $q^{\vec{l}}s_{\vec{\lambda}}$.  The vectors $\{ q^{\vec{l}}s_{\vec{\lambda}} \}_{\lambda,\vec{l}}$ form an orthonormal basis of $\mathcal{B}^r_{alg}$.\\
The map
\[
	\phi :  \mathcal{B}^r_{alg} \longrightarrow \mathcal{B}^r
\]
\[
	\phi(q^{\vec{l}}s_{\vec{\lambda}}) = [\vec{\lambda}]_{\vec{l}} \in H^{mid}_{T}(M_{\vec{l}}(r,n))
\]
 is an isometric algebra isomorphism.  This implies that the operators $H_i(n), E_i(n)$ defined above correspond under $\phi$ to multiplication by the homogeneous symmetric functions $h_i(n)$ and the elementary symmetric functions $e_i(n)$ for $n>0$ and to their adjoints for $n<0$.  Of course, since $\phi$ is an isometry and an algebra isomorphism, it is also an isomorphism of Heisenberg modules.\\
 By our construction of the Clifford algebra action on $\mathcal{A}^r$, $\phi$ is also an isomorphism of Clifford modules.  The theorem now follows from the algebraic formulation of the boson-fermion correspondence.
\end{proof}
\end{theorem}

\section{Higher Level Representations}
\subsection{A different realization of the colored bosonic Fock space}
We can repeat the constructions of Heisenberg and Clifford algebras using the space $\widetilde{\C^2/\Z_k}^{[n]}$ in place of the space $M(r,n)$, giving a different geometric realization of the the $k$-colored Fock spaces. Since the bosonic Fock space is sufficient for our purposes in this paper, we will discuss it here.\\
Recall that the fixed points
$\vec{\lambda} = (\lambda_0,\hdots,\lambda_{k-1}) \in (\widetilde{\C^2/\Z_k}^{[n]})^{\C*}$ are naturally enumerated by $k$-tuples of partitions of total size $n$.\\
For a fixed point $\vec{\lambda} = (\lambda_0,\hdots,\lambda_{k-1})$
we have a class 
$[\vec{\lambda}] = \eta'(1_{\vec{\lambda}}) \in H^{mid}_{\C^*}(\widetilde{\C^2/\Z_k}^{[n]})$, where
\[
	\eta': 	H^{mid}_{\C^*}((\widetilde{\C^2/\Z_k}^{[n]})^{\C^*} \longrightarrow 
	H^{mid}_{\C^*}(\widetilde{\C^2/\Z_k}^{[n]})
\]	
is the transport map.  Let 
$Sym^k = Sym \otimes \hdots \otimes Sym$ be the tensor product of $k$ copies of Sym, and let \\
$s_{\vec{\lambda}} = s_{\lambda_0} \otimes \hdots \otimes s_{\lambda_{k-1}}$.\\
There is an isometric algebra isomorphism
\[
	\psi : Sym^k \longrightarrow \bigoplus_n H^{mid}_{\C^*}(\widetilde{\C^2/\Z_k}^{[n]})
\]
\[
	\psi(s_{\vec{\lambda}}) = [\vec{\lambda}],
\]
where the algebra structure on  $\bigoplus_n H^{mid}_{\C^*}(\widetilde{\C^2/\Z_k}^{[n]})$ is defined exactly as it was for the space $\bigoplus_n H^{mid}_{\C^*}({\C^2}^{[n]})$.\\
For $\Sigma \subset \widetilde{\C^2/\Z_k}$ a $\C^*$ invariant curve, define the class 
$P_{\Sigma}(n) \in  \bigoplus_n H^{mid}_{\C^*}(\widetilde{\C^2/\Z_k}^{[n]})$ as the fundemental class of the variety of schemes whose support is a single point on the curve $\Sigma$.  The fundemental classes $[\Sigma]$ of $\C^*$ invariant curves $\Sigma$ span all of $H^{2}_{\C^*}(\widetilde{\C^2/\Z_k})$, so by extending linearly we can define a class $P_X(n)$ for any 
$X \in  H^{2}_{\C^*}(\widetilde{\C^2/\Z_k})$.  In particular, if $e_0,\hdots, e_{k-1}$ are the fixed points of the $\C^*$-action on $\widetilde{\C^2/\Z_k}$, we have classes $P_{\eta'(e_i)}(n) \in H^{mid}_{\C^*}(\widetilde{\C^2/\Z_k}^{[n]})$.  A detailed discussion of the $k$-dimensional Heisenberg algebra spanned by the classes $P_{\eta(e_i)}(n)$ can be found in [Q-W], which includes the following proposition.
\begin{proposition}
	The operators given by multiplication by the classes $P_{\eta(e_i)}(n)$ and their adjoints $P_{\eta'(e_i)}(-n)$ satisfy the following commutation relations:
\[
	[P_{\eta'(i)}(n),P_{\eta'(e_j)}(m)] = n \delta_{i,j} \delta_{n+m,0} Id
\]
In particular, $\bigoplus_n H^{mid}_{\C^*}(\widetilde{\C^2/\Z_k}^{[n]})$ is a geometric realization of the  $k$-colored bosonic Fock space.
\end{proposition}

\subsection{Geometric Realization of Level $k$ Representations}
The decomposition of $\coprod_{n,\vec{l}} M_{\vec{l}}(r,n)$ into connected components induces a natural grading on the vector space 
\[
	\mathcal{B}^r = \bigoplus_{n,\vec{l}} H^{mid}_{T}(M_{\vec{l}}(r,n)).
\]
For any matrix unit $e_{i,j} \in gl(r)$, the operator $e_{i,j}\otimes t^m \in \widehat{gl(r)}$ is homogeneous with respect to this geometric grading 
in the sense that if $x \in H^{mid}_{\vec{l}}(M(r,n))$ is homogeneous, then\\
 $e_{i,j}\otimes t^m(x)$ will be homogeneous as well.\\
The inclusion $\Z_k \hookrightarrow \C^* \times 1 \subset T$ induces an action of $\Z_k$ on the spaces
$M_{\vec{l}}(r,n)$ such that the connected component of the fixed point set $\coprod_{n,\vec{l}} M_{\vec{l}}(r,n)^{\Z_k}$
are $\widehat{A_{k-1}}$ quiver varieties.\\
Since the $T$- fixed points of $\coprod_{n,\vec{l}} M_{\vec{l}}(r,n)$ are the same as the $T$-fixed points of $\coprod_{n,\vec{l}} M_{\vec{l}}(r,n)^{\Z_k}$, there is a natural vector space isomorphism
\[
	\mathcal{B}^r = \bigoplus_{n,\vec{l}} H^{mid}_{T}(M_{\vec{l}}(r,n))
	\simeq \bigoplus_{n,\vec{l}} H^{mid}_{T}(M_{\vec{l}}(r,n)^{\Z_k}).
\]
Thus, the decomposition of $\coprod_{n,\vec{l}} M_{\vec{l}}(r,n)^{\Z_k}$ into connected components induces another (more refined) grading of $\mathcal{B}^r$. 
In general, elements of the form $e_{i,j} \otimes t^m$ are not homogeneous with respect to this new grading, but elements of the form $e_{i,j} \otimes t^{km}$ are homogeneous with respect to it . Thus, the corresopndences used to generate the subalgebra 
$\widehat{gl(r)}_k = gl(r) \otimes \C[t^k,t^{-k}] \oplus \C c$ can be defined naturally inside the $\Z_k$- fixed points $\coprod_{n,\vec{l}} M_{\vec{l}}(r,n)^{\Z_k} \times \coprod_{n,\vec{l}} M_{\vec{l}}(r,n)^{\Z_k}$.
Since $\widehat{gl(r)}_k$ is isomorphic to $\widehat{gl(r)}$ as Lie algebras, the geometric action of $\widehat{gl(r)}_k$ on $\bigoplus_{n,\vec{l}} H^{mid}_{T}(M_{\vec{l}}(r,n)^{\Z_k}$ will give a geometric realization of the level $k$ irreducible representations of $\widehat{gl(r)}$.\\
In the rest of this section we will use the above different realization of the $k$-colored bosonic Fock space to construct operators corresponding to the action of the homogeneous Heisenberg subalgebra $\mathfrak{h} \otimes \C[t^{k},t^{-k}] \oplus \C c \subset \widehat{gl(r)}_k$ on $\simeq \bigoplus_{n,\vec{l}} H^{mid}_{T}(M_{\vec{l}}(r,n)^{\Z_k})$.  Since the action of this subalgebra is sufficient for determining the level of the representation of the entire Lie algebra $\widehat{gl(r)}_k$, this will justify our claim that taking $\Z_k$-fixed points gives rise to level $k$ representations.

\subsection {Geometric Realization of the Heisenberg subalgebra $\mathfrak{h} \otimes \C[t^{k},t^{-k}] \oplus \C c$}
For simplicity, we first discuss the case $r=1$.\\
The fixed points $\lambda \in M(\vec{1_0},\vec{n})^{\C^*} \subset {{\C^2}^{[n]}}^{\C^*}$ are enumerated by the $k$-regular partitions, which implies that the isomorphism
\[
	\phi: Sym \longrightarrow \bigoplus_{\vec{v}} H^{mid}_{\C^*}(M(\vec{1_0},\vec{v}) \simeq
	\bigoplus_n H^{mid}_{\C^*} ({\C^2}^{[n]})
\]
restricts to an isomorphism 
\[
	\phi: Sym_{k-reg} \longrightarrow \bigoplus_{n} H^{mid}_{\C^*}M(\vec{1_0},\vec{n}).
\]
In particular, the power-sum symmetric function $p_{kn}$ is now realized geometrically as a class $P(kn)  = \phi(p_{kn})\in  H^{mid}_{\C^*} (M(\vec{1_0},\vec{n}))$.\\

We have a commutative diagram
\begin{center}
$\xymatrix{\bigoplus_nH^{mid}_{\C^*}(M(1_0,\vec{n})) \ar[r]^{f^*}  & \bigoplus_n H^{mid}_{\C^*}(\widetilde{\C^2/\Z_k}^{[n]})
\\ Sym_{k-reg} \ar[r]^{g} \ar[u]^{\phi}&Ê Sym^k \ar[u]^{\psi}}$
\end{center}

where 
\[
	f^*: H^{mid}_{\C^*}(M(1_0,\vec{n})) \longrightarrow H^{mid}_{\C^*}(\widetilde{\C^2/\Z_k}^{[n]})
\]	
 is the cohomology map induced from our $\C^*$ -equivariant diffeomorphism
 \[
 	f: \widetilde{\C^2/\Z_k}^{[n]} \longrightarrow M(1_0,\vec{n}).
 \]
 The bottom arrow $g$ takes $s_{\mu}$ to $s_{\vec{\lambda}}$ where $\vec{\lambda}$ is the $k$-quotient of the $k$-regular partition $\mu$ (see [Lei]).\\

Let $1 = \sum_{i=0}^{k-1} e_i$ be the unit in $H^0_{\C^*}(\widetilde{\C^2/\Z_k}^{\C^*})$; we have the associated Heisenberg class
\[
	P_{\eta'(1)}(n) = \sum_{i=0}^{k-1} P_{\eta'(e_i)}(n) \in H^{mid}_{\C^*}(\widetilde{\C^2/\Z_k}^{[n]}).  
\]
\begin{proposition}
	$f^*(P(kn)) = P_{\eta'(1)}(n)$
\begin{proof}
	Using the above commutative diagram, we have only to compute the image of the power-sum symmetic function $p_{kn}$ under the bottom isomorphism $g$.  This is done in [Lei], where he shows that the $g(p_{kn}) = \sum_{j=0}^{k-1} (p_n)_j$, where \\
$(p_n)_j = 1 \otimes \hdots \otimes p_n \otimes \hdots \otimes 1$ is equal to $p_n$ in the $j$th coordinate.  
It follows that 
\[
	f^*(P(kn)) = \psi \circ g \circ \phi^{-1}(P(kn)) = \psi \circ g(p_{kn}) = \psi (\sum_{j=0}^{k-1} (p_n)_j) =
	 \sum_{j=0}^{k-1} P_{\eta'(e_j)}(n) = P_{\eta'(1)}(n)
\]
\end{proof}
\end{proposition}

\begin{lemma}
Let 
$1= \sum_{i=1}^{k} e_i  \in H^0_{\C^*}(\widetilde{\C^2/\Z_k}^{\C^*})$ denote the unit, and let $\eta' : H^*_{\C^*}(\widetilde{\C^2/\Z_k}^{C^*}) \longrightarrow H^*_{\C^*}(\widetilde{\C^2/\Z_k})$ be the transport map.  Then
\[
	\langle \eta'(1), \eta'(1) \rangle = k
\]
\begin{proof}
\begin{equation*}
	\langle \eta'(1), \eta'(1) \rangle = (-1)p_* (i_*^{-1}) (\eta'(1) \cup \eta'(1)) = 
	(-1)p_* (i_*^{-1})(i_*(1) \cup i_*(1)) \cup e(N^-)^{-2} =(-1) p_*(1\cup 1 \cup e(N) \cup e(N^-)^{-2} =
\end{equation*}
\begin{equation*}
	= p_*(1 \cup 1) = p_*((\sum_{i=1}^{k} e_i) \cup (\sum_{j=1}^{k} e_j))
	= \sum_{i,j = 1}^{k} \delta_{i,j} = k
\end{equation*}
\end{proof}
\end{lemma}

\begin{corollary}
As operators on $\oplus_{\vec{v}} H^{mid}_T(M(\vec{1_0},\vec{v})$, the operators $P_{\eta'(1)}(n)$ satisfy the commutation relations
\[
	[P_{\eta'(1)}(n),P_{\eta'(1)}(m)] = kn \delta_{n+m,0} Id.
\]
\begin{proof}
Using the above lemma, we can repeat the argument of [Na], chapter 9.  See also [Vas], [Q-W].
\end{proof}
\end{corollary}
In other words, the above proposition shows that the diffeomorphism between $\widetilde{\C^2/\Z_k}^{[n]}$  and $M(1_0,\vec{n})$ identifies two different Heisenberg subalgebras:  the diagonal one-dimensional Heisenberg subalgebra spanned by the operators $ P_{\eta'(1)}(n)$ and the one-dimensional Heisenberg subalgebra $\mathfrak{h} \otimes \C[t^k,t^{-k}]$ spanned by the operators $P(kn)$.  The above corollary gives a geometric proof that these Heisenberg algebras act on the Fock space with level $k$.
\vskip0.5cm

\noindent
We now consider the case when $r \neq 1$.  Since the $\Z_k$ action on $M_{\vec{l}}(r,n)$ commutes with the $T$ action, the spaces $M_{\vec{l}}(r,n)^{\Z_k}$ are $T$-stable, and the transport map $\eta'$ gives an isomorphism
\[
	\eta': \bigoplus_{\vec{v}} H^{mid}_T(M_{\vec{l}}(r,n)^{\Z_k} \longrightarrow 
	\bigoplus_n H^{mid}_T(M_{\vec{l}}(r,n)).
\]
The $r$-dimensional subtorus $T'\subset T$ acts on $M_{\vec{l}}(r,n)^{\Z_k}$, and we have a fixed point decomposition of
$M_{\vec{l}}(r,n)^{T'\times \Z_k}$
into products of components of $({\C^2}^{[n]})^{\Z_k}$. \\
For $i = 0,\hdots, r-1$, let $P_i^k(n)$ be the operator $f^*(P_{\eta'(1)}(n))$ which acts on the $i$-th tensor product factor in
\[
	\bigoplus_{\vec{v}} H^{mid}_T(M_{\vec{l}}(r,n)^{T'\times \Z_k}) \simeq
	 \bigoplus_{\sum_{j=0}^{r-1} \vec{v_j} = \vec{v}} 
	H^{mid}_T(M(\vec{1_{l_0}},\vec{v_0}))\otimes \hdots \otimes 
	H^{mid}_T(M(\vec{1_{l_{r-1}}},\vec{v_{r-1}}))
\]

Using the isomorphism $\eta_1: H^{mid}_T(M_{\vec{l}}(r,n)^{T'\times \Z_k}) \longrightarrow H^{mid}_T(M_{\vec{l}}(r,n)^{\Z_k})$, we can transport the $P_i^k(n)$ to operators on
\[
	\bigoplus_{\vec{l},n} H^{mid}_T(M_{\vec{l}}(r,n)^{\Z_k})
\]
We will denote these operators by $P_i^k(n)$ as well.\\
In summary, on the space
\[
	\bigoplus_{\vec{l},n} H^{mid}_T(M_{\vec{l}}(r,n)^{\Z_k})
\]
we have an action of an $r$-dimensional Heisenberg algebra generated by the operators $P_i^k(n)$.  As in the case $r=1$, this Heisenberg algebra acts at level $k$.  Using the isomorphism
\[
	\bigoplus_{\vec{l},n} H^{mid}_T(M_{\vec{l}}(r,n)^{\Z_k}) \simeq \mathcal{B}^r
\]
we have the following proposition.
\begin{proposition}
	The operators $P_i^k(n)$ give a geometric realization of the action of the subalgebra $\mathfrak{h} \otimes \C[t^k,t^{-k}]$ on $\bigoplus_{\vec{l},n} H^{mid}_T(M_{\vec{l}}(r,n)^{\Z_k})$.
\begin{proof}
This follows immediately from proposition 10, as in the case $r=1$: the operators $P_i^k(n)$ correpond to multiplication by the power sum symmetric functions $p_i(kn)$ if $n > 0$ and to the adjoint operator $\frac{\partial}{\partial p_i(kn)}$ if $n < 0$.  Thus the operators $P_i(kn)$ span the subalgebra  $\mathfrak{h} \otimes \C[t^k,t^{-k}]$.
\end{proof}
\end{proposition}

\section{Level-Rank Duality}
What we have done in this paper is construct geometric realizations of Heisenberg and Clifford algebra representations, and show that if one passes to representations of $\widehat{gl(r)}$, then on the space
\[
	\mathcal{B}^r_k = \bigoplus_{n,\vec{l}} H^{mid}_T(M_{\vec{l}}(r,n)^{\Z_k})
\]
we get level $k$ representations.  On the other hand, Nakajima's original construction [Na94], [Na98]  makes the very same space $\mathcal{B}^r_k$ into a level $r$ representation of $\widehat{gl(k)}$ (strictly speaking, Nakajima constructs representations of $\widehat{sl(r)}$ of level $k$ on Borel-Moore homology, but using equivariant cohomology one can extend this to an action of the entire algebra $\widehat{gl(k)}$.)  Thus we have a geometric interpretation of a level-rank duality in the representation theory of the affine Lie algebras $\widehat{gl(n)}$.\\
Algebraically, this duality was discovered and investigated in [Fr 81], where it was applied most notably to character theory and to KdV equations .  All of the constructions in that paper should have geometric interpretations using Nakajima's construction and our dual construction.\\

The existence of a level-rank duality in Nakajima quiver varieties raises some geometric questions.  For example, instead of viewing our duality as a construction of two different representations on the space 
$\mathcal{B}^r_k$, we can observe that we now have two different geometric realizations of the same representation:  Nakajima's contruction on $\mathcal{B}^r_k$ and our construction on $\mathcal{B}^k_r$.  The author does not currently know a good geometric explanation for the exitence of two dual constructions of the same representation, and finding such an explanation seems the most important question emerging from the present paper.  In the instanton language, what we are looking for is a relationship between
\[
	U(r)\hskip0.2cm instantons \hskip0.2cm on \hskip0.2cm \widetilde{\C^2/\Z_k}
\]
and
\[
	U(k)\hskip0.2cm instantons \hskip0.2cm on \hskip0.2cm \widetilde{\C^2/\Z_r}
\]
When $k=1$, this means finding a geometric relationship between the Hilbert scheme 
$\widetilde{\C^2/\Z_r}^{[n]}$ and the space $M(r,n)$.  Note that both spaces have torus actions with isolated fixed points, and in both cases the fixed point sets are parametrized by $r$-tuples of partitions of total size $n$.

We should also remark that our decision to restrict our discussion to the type $\widehat{A}$ case is not merely for convenience, despite the fact that Nakajima's construction works for more arbitrary $\Gamma \subset SL(2,\C)$.  Indeed, it is tempting to try to compare the instanton moduli spaces
\[
	G_1\hskip0.2cm instantons \hskip0.2cm on \hskip0.2cm \widetilde{\C^2/\Gamma_2}
\]
\[
	G_2\hskip0.2cm instantons \hskip0.2cm on \hskip0.2cm \widetilde{\C^2/\Gamma_1}
\]
where  $G_i$ is related to $\Gamma_i$, $i=1,2$ by the McKay correspondence.
Unfortunately these instanton moduli spaces are not complete, though in the case where $G$ is of type $A$ we can use the completion $M(r,n)$, known as the "Geiseker compactification".  If $G$ is not of type $A$, the author is not aware of a smooth completion of the moduli space of $G$-instantons on a surface.
There is however, the singular "Uhlenbeck compactification", which does exist for more general $G$.  One is tempted to conjecture that there should be a "level-rank" duality between the intersection cohomologies of the appropriate Uhlenbeck compactifications.

\section{Bibliography}
[C-K]		D.A. Cox, S. Katz, Mirror Symmetry and Algebraic Geometry, Mathematical Surveys and Monographs 68, AMS (1999)\\\\
\noindent
[Fr-K]		I. Frenkel, V. Kac	Basic representations of affine Lie algebras and dual resonance  models.  Invent. Math. 62 (1980/81), no. 1, 23--66. \\\\
\noindent
[Fr 82]		I. Frenkel, Representations of affine Lie algebras, Hecke modular forms and Kerteweg-De Vries type equations, in Lie Algebras and Related Topics, Lecture Notes in Mathematics, no. 933,
D. Winter, ed. (Springer-Verlag, Berlin, 1982)\\\\
\noindent
[Fr 81]		 I. Frenkel, 	Two constructions of affine Lie algebra representations and boson-fermion corre-spondence in quantum field theory, J. Funct. Anal. 44 (1981) 259-327.\\\\
\noindent
[Fr-Sa]	I. Frenkel, A. Savage,	Bases of representations of type A affine Lie algebras via quiver varieties and statistical mechanics, Inter. Math. Res. Notices 28 (2003), 1521-1547\\\\
\noindent
[GK]		V.Ginzburg, M.Kapranov, 	Hilbert Schemes and Nakajima's Quiver varieties, unpublished manuscript.\\\\
\noindent
[Gr]		I. Grojnowski,	Instantons and affine Lie algebras I: the Hilbert Scheme and vertex operators, Math. Res. Letters 3 (1996), 275-291\\\\
\noindent
[Hai]		M. Haiman, 	Combinatorics, symmetric functions and Hilbert schemes
In CDM 2002: Current Developments in Mathematics in Honor of  Wilfried Schmid and George Lusztig, International Press Books  (2003) 39-112.\\\\
\noindent
[It-Na]		Ito, Nakamura,	McKay Correspondence and Hilbert Schemes, Proc. Jap. Acad. 72 (1996)135-138\\\\
\noindent
[K-K-L-W]		Kac, Kazhdan, Lepowsky, Wilson	Realization of the basic representations of the Euclidean Lie  algebras.  Adv. in Math. 42 (1981), no. 1, 83--112. \\\\
\noindent
[Kuz]		A. Kuznetsov . Quiver varieties and Hilbert schemes. arXiv:math.AG/0111092.\\\\
\noindent
[Lei]		S. Leidwanger, Basic representations of  and  and the combinatorics of partitions, Adv. in Math. 141 (1999), 119-154.\\\\
\noindent
[Le-Le]		B. Leclerc, S. Leidwanger, Schur functions and affine Lie algebras, J. Algebra 210 (1998), 103-144.\\\\
\noindent
[Mac]	Macdonald I.G., Symmetric functions and Hall polynomials, 2nd edition, Oxford University Press, Oxford, 1995. 18\\\\
\noindent
[Na 94]	H. Nakajima, 	Instantons on ALE spaces, quiver varieties, and Kac-Moody  algebras, Duke Math. 76 (1994) 365--416.\\\\
\noindent
[Na 98]	H. Nakajima, 	Quiver Varieties and Kac-Moody Algebras, Duke Math., 91, (1998), 515--560.\\\\
\noindent
[Na 97]	H. Nakajima,	Heisenberg Algebra and Hilbert Schemes of Points on Projective Surfaces, Ann. of Math. 145, (1997) 379--388.\\\\
\noindent
[Na 99]	H. Nakajima,	Lectures on Hilbert schemes of points on surfaces, AMS Univ. Lecture Series, 1999\\\\
\noindent
[Na]		H. Nakajima,	Jack polynomials and Hilbert schemes of points on surfaces, unpublished preprint,  alg-geom/9610021.\\\\
 \noindent
[NY1]	Instanton counting on blowup. I. $4$-dimensional pure gauge theory (with Kota Yoshioka), Invent. Math 162 (2005), no. 2, 313--355;\\\\
\noindent
[NY2]	Instanton counting on blowup. II. $K$-theoretic partition funtion (with Kota Yoshioka), Transform. Groups 10 (2005), no. 3-4, 489--519;\\\\
\noindent
[PS]		Pressley, Segal	Loop Groups	 Oxford Mathematical Monographs. Oxford Science Publications. The Clarendon Press, Oxford University Press, New York, 1986.\\\\
\noindent
[Sa]		A geometric boson-fermion correspondence, math.RT/0508438\\\\
\noindent
[Se]		Segal 	Unitary representations of some infinite-dimensional groups. Comm. Math. Phys. 80 (1981), no. 3, 301--342.\\\\
\noindent
[Vas]		E. Vasserot	 Sur lÕanneau de cohomologie du schema de Hilbert de C2
, C.-R. Acad. Sc. Paris 332 (2001), 7Ð12.\\\\
\noindent
[Q-W]		Qin, W. Wang	Hilbert Schemes of Points on the Minimal Resolution and Soliton Equations,
math.QA/0404540\\\\
\noindent

contact:	Tony Licata\\
anthony.licata@yale.edu

\end{document}